\documentclass[letterpaper,11pt]{amsart}

\usepackage[all]{xy}                        %

\CompileMatrices                            

\UseTips                                    

\usepackage[bookmarks=true]{hyperref}       

\usepackage{amssymb,latexsym,amsmath,amscd}
\usepackage{xspace}

\usepackage{graphicx}

\reversemarginpar

\theoremstyle{plain}
\newtheorem{theorem}{Theorem}[section]
\newtheorem*{theorem*}{Theorem}
\newtheorem{proposition}[theorem]{Proposition}
\newtheorem{corollary}[theorem]{Corollary}
\newtheorem{lemma}[theorem]{Lemma}

\theoremstyle{definition}
\newtheorem{definition}[theorem]{Definition}
\newtheorem{exercise}[theorem]{Exercise}

\newtheorem{remark}[theorem]{Remark}
\newtheorem{example}[theorem]{Example}

\newcommand{\enm}[1]{\ensuremath{#1}}          %
\newcommand{\op}[1]{\operatorname{#1}}
\newcommand{\cal}[1]{\mathcal{#1}}

\renewcommand{\bar}[1]{\overline{#1}}

\newcommand{\CC}{\enm{\mathbb{C}}}

\newcommand{\NN}{\enm{\mathbb{N}}}
\newcommand{\RR}{\enm{\mathbb{R}}}
\newcommand{\QQ}{\enm{\mathbb{Q}}}
\newcommand{\ZZ}{\enm{\mathbb{Z}}}

\renewcommand{\AA}{\enm{\mathbb{A}}}
\newcommand{\PP}{\enm{\mathbb{P}}}

\newcommand{\Jj}{\enm{\cal{J}}}

\newcommand{\Oo}{\enm{\cal{O}}}

\renewcommand{\phi}{\varphi}
\renewcommand{\theta}{\vartheta}
\renewcommand{\epsilon}{\varepsilon}

\newcommand{\Spec}{\op{Spec}}

\newcommand{\Supp}{\op{Supp}}
\newcommand{\codim}{\op{codim}}

\newcommand{\tensor}{\otimes}         

\newcommand{\abs}[1]{\left\vert#1\right\vert}
\newcommand{\set}[1]{\left\{#1\right\}}

\renewcommand{\to}[1][]{\xrightarrow{\ #1\ }}

\newcommand{\defeq}{\stackrel{\scriptscriptstyle
\op{def}}{=}}

\newcommand{\ie}{\textit{i.e.},\ }           

\newcommand{\Bl}{\mathit{Bl}}
\renewcommand{\aa}{\mathfrak{a}}
\newcommand{\bb}{\mathfrak{b}}
\newcommand{\mm}{\mathfrak{m}}
\newcommand{\ord}{\op{ord}}
\renewcommand{\div}{\op{div}}
\newcommand{\floor}[1]{\lfloor #1 \rfloor}
\newcommand{\ceil}[1]{\lceil #1 \rceil}
\newcommand{\lct}{\op{lct}}
\newcommand{\Newt}{\op{Newt}}
\newcommand{\old}[1]{}

\newcommand{\MI}[1]{\Jj ({#1})}
\newcommand{\fra}{\aa}
\newcommand{\frb}{\bb}
\newcommand{\frq}{\mathfrak{q}}

\newcommand{\lra}{\longrightarrow}
\newcommand{\noi}{\noindent}
\newcommand{\bull}{_{\bullet}}

\begin{document}

\title[An Informal Introduction to Multiplier Ideals]{An Informal
introduction to
\\ multiplier ideals}
\author{Manuel Blickle}
\address{Universit\"at Essen, FB6 Mathematik, 45117
Essen, Germany}
\email{manuel.blickle@uni-essen.de}
\author{Robert Lazarsfeld}
\address{University of Michigan, Ann Arbor, MI 48109, USA}
\email{rlaz@umich.edu}
\thanks{Research of the second author partially supported by NSF Grant DMS
0139713}

\maketitle

\section{Introduction}

Given a smooth complex variety $X$  and an ideal (or
ideal sheaf) $\fra$ on $X$, one can attach to $\fra$ a
collection of \textit{multiplier ideals} $\MI{\fra^c}$
depending on a rational weighting parameter $c > 0$.
These ideals, and the vanishing theorems they satisfy,
have found many applications in recent years. In the
global setting they have been used to study
pluricanonical and other linear series on a projective
variety (\cite{Demailly93c}, \cite{Angehrn-Siu95a}, \cite{Siu98a},
\cite{Ein-Lazarsfeld97a}, \cite{ELNull}, \cite{Demailly99b}). More recently
they have led to the discovery of some surprising uniform results in local
algebra (\cite{ELS1},
\cite{ELS2}, \cite{ELSV}). The purpose of these lectures is to
give an easy-going and gentle introduction to the
algebraically-oriented local side of the theory.

Multiplier ideals can be approached  (and historically emerged) from three
different viewpoints. In commutative algebra they were introduced and studied
by Lipman \cite{lip.adj} in connection with the Brian\c con-Skoda
theorem.\footnote{Lipman used the term ``adjoint ideals", but this has come to
refer to a different construction.} On the analytic side of the field, Nadel
\cite{Nadel90} attached a multiplier ideal to any plurisubharmonic function,
and proved a Kodaira-type vanishing theorem for them.\footnote{In fact, the
``multiplier"  in the name refers to their analytic construction (see \S
\ref{Analytic.Interp.Subsectn})} This machine was developed and applied
with great success by Demailly, Siu and others. Algebro-geometrically, the
foundations were laid in passing by Esnault and Viehweg  in connection with
their work involving the Kawamata-Viehweg vanishing theorem. More systematic
developments of the geometric theory were subsequently undertaken by Ein,
Kawamata and the second author. We will take the geometric approach here.

The present notes   follow  closely  a short course
on multiplier ideals given by the second author at the
Introductory Workshop for the Commutative Algebra
Program at the MSRI in September
2002\footnote{Handwritten notes and the lectures on
streaming video are available at \\
http://www.msri.org/publications/video/index05.html}.
The three main lectures were supplemented with
a presentation by the first author on multiplier ideals
associated to monomial ideals (which appears here in \S
3).   We have tried to preserve in this write-up the
informal tone of these talks: thus we emphasize
simplicity over generality in statements of results,
and we present very few proofs. Our primary hope  is
to give the reader  a feeling for what multiplier ideals
are and how they are used. For a detailed development
of the theory from an algebro-geometric perspective we
refer to Part Three of the forthcoming book
\cite{PAG}. The analytic
picture is covered in Demailly's lectures
\cite{Dem.Mult}.

We conclude this Introduction by fixing the set-up in
which we work and giving a brief preview of what is to
come.
 Throughout these notes, $X$ denotes a smooth
affine variety over an algebraically closed field $k$ of characteristic zero
and $R = k[X]$ is the coordinate ring of $X$, so that $X = \Spec R$.  We
consider a non-zero ideal $\fra \subseteq k[X]$ (or equivalently a sheaf of
ideals $\aa \subseteq \Oo_X$). Given  a rational number $c \geq 0$ our plan is
to define and study the multiplier ideal
\[
    \Jj(c \cdot \aa)\ =\ \Jj(\aa^c) \ \subseteq \ k[X].
\] As we proceed, there are two ideas to keep in
mind. The first is that $\Jj(\aa^c)$ measures in a somewhat subtle manner the
singularities of the divisor of a typical function $f$ in $\aa$: for fixed $c$,
``nastier" singularities are reflected by
  ``deeper" multiplier ideals. Secondly, $\MI{\fra^c}$
enjoys remarkable formal properties
   arising from the
Kawamata-Viehweg-Nadel Vanishing theorem. One can view the power of multiplier
ideals as arising from the confluence of these facts.

The theory of multiplier ideals described here has striking parallels with the
theory of tight closure developed by Hochster and Huneke in positive
characteristic. Many of the uniform local results that can be established
geometrically via multiplier ideals can also be proven (in more general
algebraic settings) via tight closure. For some time the actual connections
between the two theories were not well understood. However very recent work
\cite{HaraYosh}, \cite{Takagi.MultTest} of Hara-Yoshida and Takagi  has
generalized tight closure theory to define a so called test ideal $\tau(\aa)$,
which corresponds to the multiplier ideal $\Jj(\aa)$ under reduction to
positive characteristic. This provides a first big step towards identifying
concretely the links between these theories.

Concerning the organization of these notes, we start in \S 2 by giving the
basic definition and several examples. Multiplier ideals of monomial ideals are
discussed in detail in \S 3. Invariants arising from multiplier ideals, with
some applications to uniform Artin-Rees numbers, are taken up in \S 4. Section
5 is devoted to a discussion of some basic  results about multiplier ideals,
notably Skoda's theorem and the restriction and subaddivity theorems. We
consider asymptotic constructions in \S 6, with applications to uniform bounds
for symbolic powers following \cite{ELS1}.

We are grateful to Karen Smith for suggestions concerning these notes.

\section{Definition and Examples}\label{sec.defex} As
just stated, $X$ is a smooth affine variety   of dimension $n$ over an
algebraically closed field of characteristic zero, and we fix an ideal $\fra
\subseteq k[X]$ in the coordinate ring of $X$.  Very little is lost by focusing
on the case $X = \CC^n$ of affine $n$-space over the complex numbers $\CC$, so
that $\fra \subseteq \CC[x_1, \ldots, x_n]$ is an ideal in the polynomial ring
in $n$ variables.

\subsection{Log resolution of an ideal} The starting
point is to realize the ideal $\fra$ geometrically.

\begin{definition} A \emph{log resolution} of an ideal
sheaf $\aa
\subseteq \Oo_X$ is a proper, birational map $\mu: Y
\to X$ whose exceptional locus is a  divisor $E = \text{Exceptional}(\mu)$
such that
\begin{enumerate}
    \item $Y$ is non-singular.
    \item $\aa \cdot \Oo_{Y} =\mu^{-1}\aa= \Oo_{Y}(-F)$
with $F=\sum r_iE_i$ an
    effective divisor.
    \item $F+E$ has simple normal crossing support.
\end{enumerate}
\end{definition}
\noi  Recall that a (Weil) divisor $D=\sum \alpha_i D_i$ has simple normal
crossing support if each of its irreducible components $D_i$ is smooth, and if
locally analytically one has coordinates $x_1,\ldots,x_n$ of $Y$ such that
$\Supp D=\sum D_i$ is defined by $x_1\cdot\ldots\cdot x_a$ for some $a$ between
$1$ and $n$. In other words, all the irreducible components of $D$ are smooth
and intersect transversally. The existence of a log resolution for any sheaf of
ideals in any variety over a field of characteristic zero is essentially
Hironaka's celebrated result on resolution of singularities
\cite{Hironaka.ResSing}. Nowadays there are more elementary constructions of
such resolutions, for instance \cite{Bierstone-Milman97},
\cite{EncVill.Desing} or \cite{Paranjape}.

\begin{figure}[h]
\vskip 5pt
\includegraphics{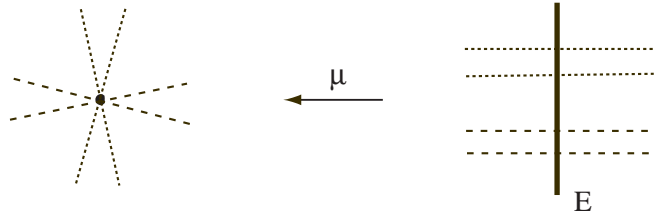}
 \vskip 5pt
\caption{Log resolution of $(x^2, y^2)$} \label{Resolve.Double.Pt}
\end{figure}

\begin{example} Let $X=\AA^2=\Spec k[x,y]$ and $\aa =
(x^2,y^2)$. Blowing up the origin in $\AA^2$ yields
\[
    Y = \Bl_0(\AA^2) \to[\mu] \AA^2=X.
\] Clearly, $Y$ is nonsingular. Computing on the chart
for which the blowup $\mu$ is a map from $\AA^2 \to \AA^2$ given by $(u,v)
\mapsto (u,uv)$ shows that $\aa \cdot \Oo_{Y} = \Oo_{Y}(-2E)$. On the described
chart we have $\aa \cdot \Oo_Y = (u^2,u^2v^2)=(u^2)$ and $(u=0)$ is the
equation of the exceptional divisor. This resolution is illustrated in Figure
\ref{Resolve.Double.Pt}, where we have  drawn schematically the curves in
$\AA^2$ defined by typical $k$-linear combinations of generators of $\fra$, and
the proper transforms of these curves on $Y$. Note that these proper transforms
do not meet: this reflects the fact that $\fra$ has become principal on
$Y$.
\end{example}

\begin{figure}[h]
\vskip 5 pt
\includegraphics{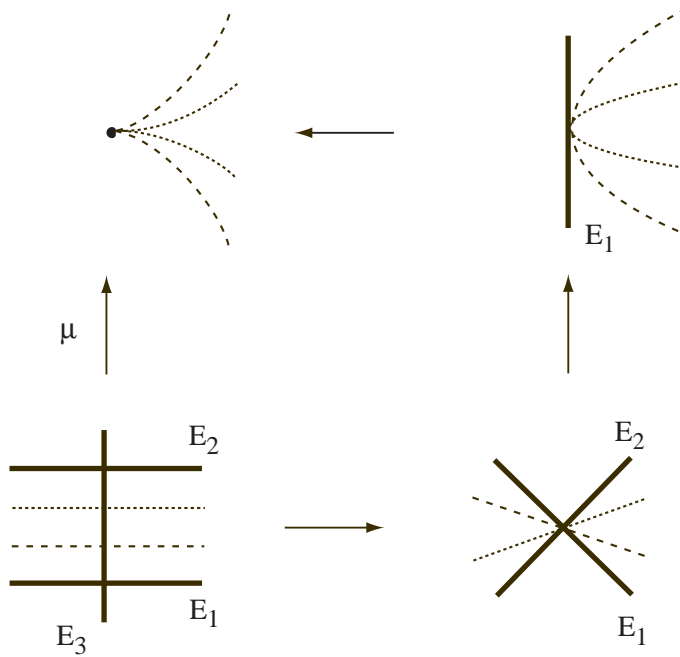}
\caption{Log resolution of $(x^3,y^2)$.} \label{Resolve.Cusp} \vskip 10 pt
\end{figure}

\begin{example} Now let $\aa = (x^3,y^2)$. In this case
a log resolution is constructed by the familiar
sequence of three blowings-up used to resolve a
cuspidal curve (Figure \ref{Resolve.Cusp}).
 Here we have $\aa \cdot \Oo_{Y} =
\Oo_{Y}(-2E_1-3E_2-6E_3)$ where $E_i$ is the
exceptional divisor of the
$i$th blowup.
\end{example}

These examples illustrate the principle that a  log
resolution of an ideal
$\aa$ is very close to being the same as a resolution
of singularities of a divisor of a general function in
$\aa$.

\subsection{Definition of multiplier ideals} Besides a
log resolution of $\mu:Y  \to X$ of the ideal $\aa$, the other ingredient for
defining the multiplier ideal is the relative canonical divisor
\[
    K_{Y/X}=K_{Y}-\mu^*K_X= \div  (\det (\op{Jac} \mu)).
\] It is unique as a divisor (and not just as a divisor
class) if one requires its support to be contained in the exceptional locus of
$\mu$. Alternatively, $K_{Y/X}$ is the effective divisor defined by the
vanishing of the determinant of the Jacobian of $\mu$. The canonical divisor
$K_X$ is just the class corresponding to the canonical line bundle $\omega_X$.
If $X$ is smooth, $\omega_X$ is just the sheaf of top differential forms
$\Omega^n_X$ on $X$.

Extremely useful for basic computations of multiplier
ideals is the following proposition, see \cite{Ha},
Exercise II.8.5.

\begin{proposition}\label{prop.CanBlow} Let $Y=\Bl_Z X$
where $Z$ is a smooth subvariety of the smooth variety
$X$ of codimension $c$. Then the relative canonical
divisor $K_{Y/X}$ is $(c-1)E$,   $E$ being the
exceptional divisor of the blowup.
\end{proposition} Now we can give a provisional
definition of the multiplier ideal of an ideal $\fra$:
it coincides in our setting with Lipman's construction
in
\cite{lip.adj}.
\begin{definition} Let $\aa \subseteq k[X]$ be an
ideal.  Fix a log resolution $\mu: Y \to X$ of $\aa$ such that
$\aa\cdot\Oo_{Y}=\Oo_{Y}(-F)$, where $F=\sum r_iE_i$, and $K_{Y/X}=\sum
b_iE_i$. The \emph{multiplier ideal} of $\aa$ is
\begin{equation*}
\begin{split}
    \Jj(\aa) \ &=\ \mu_*\Oo_{Y}(K_{Y/X} - F) \\
              &= \ \big \{ \,  h \in k[X] \mid
\div(\mu^*h)+K_{Y/X} - F \, \geq \, 0 \, \big \}
             \\
             &=\ \big\{\, h \in k[X] \mid \ord_{E_i}(\mu^* h)
\geq r_i - b_i
             \text{ for all $i$}\, \big\}.
\end{split}
\end{equation*}
(We will observe later that this is independent of the
choice of resolution.)
\end{definition}

The definition may seem at first blush a little mysterious. One way to motivate
it is to note that $\MI{\fra}$ is the push-forward of a bundle which is very
natural from the viewpoint of vanishing theorems. In fact, the bundle
$\Oo_Y(-F)$ appearing above is (close to being) ample for the map $\mu$.
Therefore $K_{Y/X} - F$ has the shape to which Kodaira-type vanishing results
will apply. In any event, the definition will justify itself before long
through the properties of the ideals so defined.

\begin{exercise} Use the fact that
$\mu_*\omega_Y=\omega_X$ to show that
$\Jj(\aa)$ is indeed an ideal in $k[X]$.
\end{exercise}
\begin{exercise}\label{ex.intclo} Show that the
integral closure $\bar{\aa}$ of $\aa$ is equal to $\mu_*
\Oo_{Y}(-F)$. Use this to conclude that $\aa \subseteq
\bar{\aa} \subseteq
\Jj(\aa) = \bar{\Jj(\aa)}$. (Recall that the integral
closure of an ideal
$\aa$ consists of all elements $f$ such that $v(f) \geq
v(a)$ for all valuations $v$ of $\Oo_X$.)
\end{exercise}

\begin{exercise}  \label{Mult.Ideal.Int.Closure.Ideal}
Verify that for ideals $\aa \subseteq \bb$ one has $\Jj(\aa) \subseteq
\Jj(\bb)$. Use this and the previous exercise to show that
$\Jj(\aa)=\Jj(\bar{\aa})$.
\end{exercise}

The above definition of the multiplier ideal is not
 general enough for the most interesting applications.
As it turns out, allowing an additional
rational (or real) parameter
$c$   considerably increases the
power of the theory.

Note that a log resolution of an
ideal
$\aa$ is at the same time a log resolution of any
integer power $\aa^n$ of that ideal. Thus we extend the
last definition, using the same log resolution for
every $c \geq 0$:
\begin{definition}\label{Gen.Mult.Ideal.Def}
For every rational number $c\geq 0$, the \emph{multiplier ideal} of the ideal
$\aa$ with exponent (or coefficient) $c$ is
\begin{equation}
\begin{split}
    \Jj(\aa^c)\ = \ \Jj(c\cdot \aa) \ &= \
\mu_*\Oo_{Y}(K_{Y/X}-\floor{c\cdot
    F}) \\
    &= \ \big \{ h \in k[X] \, |\, \ord_{E_i}(\mu^*h) \,
\geq \,
\floor{cr_i}-b_i \text{
    for all $i$}\, \big \}
\end{split}
\end{equation} where $\mu : Y \to X$ is a log
resolution of $\aa$ such that $\aa \cdot \Oo_Y = \Oo_Y(-F)$.
\end{definition} Note that we do not assign any meaning
to $\aa^c$ itself, only to $\Jj(\aa^c)$.\footnote{There is a way to define the
integral closure of an ideal $\aa^c$, for $c \geq 0$ rational, such that it is
consistent with the definition of the multiplier ideal. For $c=p/q$ with
positive integers $p$ and $q$, set $f \in \aa^{p/q}$ if and only if $f^q \in
\bar{\aa^p}$, where the bar denotes the integral closure.} The round-down
operation $\floor{\cdot}$ applied to a $\QQ$-divisor $D = \sum a_iD_i$ for distinct
prime divisors $D_i$ is just rounding down the coefficients. That is,
$\floor{D} = \sum \floor{a_i}D_i$. The round up $\ceil D=-\floor{-D}$ is
defined analogously.
\begin{exercise}[Caution with rounding]
    Show that rounding  does not in general commute with
restriction or pullback.
\end{exercise}

\begin{exercise} \label{MI.Max.Idea.Ex}
    Let $\mm$ be the maximal ideal of a point $x \in X$. Show that
    \[
        \Jj(\mm^c)=\begin{cases}
                        \mm^{\floor{c}+1-n} & \text{for
$c\geq n=\dim X$.} \\
                        \Oo_X               &
\text{otherwise.}
                   \end{cases}
    \]
\end{exercise}

\begin{example}
    Let $\aa = (x^2,y^2) \subseteq k[x,y]$. For the log
resolution of
$\aa$ as
    calculated above we have $K_{Y/X}=E$. Therefore,
    \[
        \Jj(\aa^c) = \mu_*(\Oo_{Y}(E-\floor{2c}E))=
        (x,y)^{\floor{2c}-1}
    \]
(In view of Exercise
\ref{Mult.Ideal.Int.Closure.Ideal}, this  is
 a special case of Exercise \ref{MI.Max.Idea.Ex}.)
\end{example}
\begin{example}
    Let $\aa = (x^2,y^3)$. In this case we computed a
log
    resolution with $F= 2E_1+3E_2+6E_3$. Using the
basic formula
    for the relative canonical divisor of a blowup
along a smooth
    center, one computes $K_{Y/X}= E_1 + 2 E_2 +4 E_3$.
    Therefore,
    \[
    \begin{split}
        \Jj(\aa^c)  &=\ \mu_*\big(\Oo_{Y}(E_1 + 2 E_2 +4
        E_3-\floor{c(2E_1+3E_2+6E_3)})) \\
                   &=
\mu_*(\Oo_{Y}((1-\floor{2c})E_1+(2-\floor{3c})E_2+(4-\floor{6c})E_3)).
    \end{split}
    \]
    This computation shows that for $c < 5/6$ the
multiplier ideal
    is trivial, \ie $\Jj(\aa^c)=\Oo_X$. Furthermore,
    $\Jj(\aa^{\frac{5}{6}})=(x,y)$. The next
coefficient for
    which the multiplier ideal changes is $c=1$. This
behavior of
    multiplier ideals to be piecewise constant with
discrete jumps
    is true in general and will be discussed in more
detail later.
\end{example}

\begin{exercise} [Smooth ideals]
\label{Smooth.Ideal.Exercise} Suppose that
$\frq
\subseteq k[X]$ is the ideal of a smooth subvariety $Z
\subseteq X$ of pure codimension $e$. Then
\[ \MI{\frq^\ell} \ = \ \frq^{\ell + 1 - e}. \]
(Blowing up $X$ along $Z$ yields a log resolution of
$\frq$.) The case of fractional exponents is similar.
\end{exercise}

\subsection{Two basic properties}
The definitions of the previous subsection
are   justified by the fact that they lead
to two very basic results.

The first point is that the ideal $\MI{\fra^c}$
constructed in Definition \ref{Gen.Mult.Ideal.Def} is
actually independent of the choice of resolution.
\begin{theorem}
   If
    $X_1 \to[\mu_1] X$ and $X_2 \to[\mu_2] X$  are two
log resolutions of the ideal $\aa \subseteq
\Oo_X$ such that
    $\aa \Oo_{X_i} = \Oo_{X_i}(-F_i)$, then
\[ {\mu_1}_*\big(\Oo_{X_1}(K_{X_1/X}-\floor{c\cdot
    F_1}\big)\ =
\ {\mu_2}_*\big(\Oo_{X_2}(K_{X_2/X}-\floor{c\cdot
    F_2}\big).\]
\end{theorem}
\noi  As one would expect, the proof involves dominating $\mu_1$ and $\mu_2$ by
a third resolution. It is in the course of this argument that it becomes
important to know that $F_1$ and $F_2$ have normal crossing support, see
\cite[Chapter 9]{PAG}.

\begin{exercise}
    By contrast, give an example to show that if $c$
is non-integral, then the ideal
$\mu_*(-\floor{cF})$ may
    indeed depend on the  log resolution $\mu$.
\end{exercise}

The second fundamental fact is a vanishing theorem for
the sheaves computing multiplier ideals.
\begin{theorem} [Local Vanishing Theorem]
\label{Local.Vanishing.Thm}
    Consider an ideal $\fra \subseteq k[X]$ as above,
and let
$\mu: Y
\to X$ be a log resolution of $\aa$ with
    $\aa\cdot \Oo_Y = \Oo_Y(-F)$. Then
    \[
        R^i\mu_*\Oo_Y(K_{Y/X}-\floor{cF}) = 0
    \]
    for all $i > 0$ and $c > 0$.
\end{theorem} \noi This leads one to expect that the
multiplier ideal, being the zeroth derived image of
$\Oo_Y(K_{Y/X}-\floor{cF})$ under $\mu_*$, will display
particularly good cohomological properties.

Theorem \ref{Local.Vanishing.Thm} is a special case of the Kawamata-Viehweg
vanishing theorem for a mapping, see \cite[Chapter 9]{PAG} . It is the
essential fact underlying all the applications of multiplier ideals appearing
in these notes. When  $c$ is a natural number, the result can be seen as a
slight generalization of the classical Grauert-Riemenschneider Vanishing
Theorem. However as we shall see it is precisely the possibility of working
with non-integral $c$ that  opens the door to  applications of a non-classical
nature.

\subsection{Analytic construction of multiplier ideals}
\label{Analytic.Interp.Subsectn} We sketch briefly the
analytic construction of multiplier ideals. Let
$X$ be a smooth complex affine variety,
and $\fra \subseteq \CC[X]$ an ideal. Choose generators
$g_1,\ldots,g_p \in \fra$. Then
\[
    \Jj(\aa^c)^{\text{an}}\  =_{\text{locally}} \ \Big
\{ \, h \text{ holomorphic }\Big | \,\,
    \frac{|h|^2}{\big(\sum |g_i|^{2}\big)^c} \text{ is
locally integrable} \, \Big
    \}.
\] In other words, the analytic ideal associated to $\MI{\fra^c}$ arises
as a sheaf of ``multipliers". See \cite[(5.9)]{Demailly99b} or
\cite[Chapter 9.3.D]{PAG} for the proof. In brief the idea is to show that
both the algebraic and the analytic definitions lead to ideals that
transform the same way under birational maps. This reduces one to the
situation where
$\fra$ is the principal ideal generated by a single monomial in local
coordinates. Here the stated equality can be checked by an explicit
calculation.

\subsection{Multiplier ideals via tight closure}
As already hinted at in the introduction there is an intriguing parallel
between effective results in local algebra obtained via multiplier ideals on
the one hand and tight closure methods on the other. Almost all the results we
will discuss in these notes are of this kind: there are tight closure versions
of the Brian\c{c}on-Skoda theorem, the uniform Artin-Rees lemma and even of the
result on symbolic powers we present as an application of the asymptotic
multiplier ideals in Section \ref{sec.symb}.\footnote{The tight closure
analogues of these result can be found in \cite{HH.BrianSkoda},
\cite{Huneke.UniformBounds} and \cite{HH.ComOrdSymbPow}, respectively.} There
is little understanding for why such different techniques (characteristic zero,
analytic in origin vs.\
 positive characteristic) seem to be tailor made to prove the
same results.

Recently, Hara-Yoshida and Takagi strengthened this parallel by
constructing in
\cite{HaraYosh} and \cite{Takagi.MultTest,TakWat,HaTak,Tak.Inv}
multiplier-like ideals using techniques modelled after tight closure theory.
Their construction builds on earlier work of Smith \cite{Smith.MultTest} and
Hara \cite{Hara.GeomIntTest}, who had established a connection between the
multiplier ideal associated to the unit ideal $(1)$ on certain singular
varieties with the so-called test ideal in tight closure. The setting of the
work of Hara and Yoshida is a regular\footnote{One feature of their theory is
that there is no reference to resolutions of singularities. As a consequence no
restriction on the singularity of $R$ arises, whereas for multiplier ideals at
least some sort of $\QQ$--Gorenstein assumption is needed.} local ring $R$ of
positive characteristic $p$. For simplicity one might again assume $R$ is the
local ring of a point in $\AA^n$. Just as with multiplier ideals, one assigns
to an ideal $\aa \subseteq R$ and a rational parameter $c \geq 0$, the
\emph{test ideal}
\[
   \tau(\aa^c) = \{\, h \in R\,|\, h I^{*\aa^c} \subseteq I \text{ for all ideals }
   I\,\}.
\]
Here $I^{*\aa^c}$ denotes the $\aa^c$--tight closure of an ideal, specifically
introduced for the purpose of constructing these test ideals
$\tau(\aa^c)$.\footnote{Similarly as for tight closure, $x \in I^{*\aa^c}$ if
there is a $h \neq 0$ such that for all $q=p^e$ one has $hx^q\aa^{\ceil{qc}}
\subseteq I^{[q]}$. Note that $I^{[q]}$ denotes the ideal generated by all
$q$th powers of the elements of $I$, whereas $\aa^{\ceil{qc}}$ is the usual
$\ceil{qc}$th power of $\aa$.} The properties the test ideals enjoy are
strikingly similar to those of the multiplier ideal in characteristic zero: For
example the Restriction Theorem (Theorem \ref{thm.restr}) and Subadditivity
(Theorem \ref{thm.subadd}) hold. What makes the test ideal a true analog of the
multiplier ideal is that under the process of reduction to positive
characteristic the multiplier ideal $\Jj(\aa^c)$ corresponds to the test ideal
$\tau(\aa^c)$, or more precisely to the test ideal of the reduction mod $p$ of
$\aa^c$.

\section{The multiplier ideal of monomial ideals}\label{sec.nonomial}
 Even though multiplier
ideals enjoy extremely good formal properties, they are very hard to compute in
general. An important exception is the class of monomial ideals, whose
multiplier ideals are described by a simple combinatorial formula, established
by  Howald \cite{Howald.MultMon}. By way of illustration we discuss this result
in detail.

To state the result let $\aa \subseteq k[x_1,\ldots,x_n]$ be a nonomial ideal,
that is an ideal generated by monomials of the form
$x^m=x_1^{m_1}\cdot\ldots\cdot x_n^{m_n}$ for $m \in \ZZ^n \subseteq \RR^n$. In
this way we can identify a monomial ideal $\aa$ of $k[x_1,\ldots,x_n]$ with the
set of exponents (contained in $\ZZ^n$) of the monomials in $\aa$. The convex
hull of this set in $\RR^n=\ZZ^n \tensor \RR$ is called the \emph{Newton
polytope} of $\aa$ and it is denoted by $\Newt(\aa)$. Now Howald's result
states:
\begin{theorem}\label{thm.mon}
    Let $\aa \subseteq k[x_1,\ldots,x_n]$ be a monomial
ideal. Then for every
    $c>0$,
    \[
        \Jj(\aa^c) = \langle x^m \,|\, m+(1,\ldots,1)
\in
        \text{ interior of }c \cdot \Newt(\aa) \rangle
    \]
\end{theorem}

\begin{figure}[h]\label{fig.NPx4}
\includegraphics{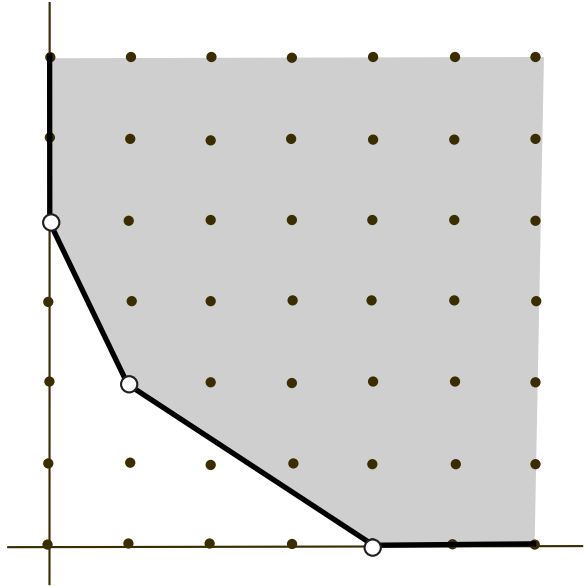}
\caption{Newton polytope of $(x^4,xy^2,y^4)$ }
\end{figure}

For example, the picture of the Newton polytope of the monomial ideal $\aa =
(x^4,xy^2,y^4)$ in Figure \ref{fig.NPx4} shows, using Howald's result,
that
$\Jj(\aa)=(x^2,xy,y^2)$. Note that even though $(0,1)+(1,1)$ lies in the Newton
polytope $\Newt(\aa)$ it does not lie in the interior. Therefore, the monomial
$y$ corresponding to $(0,1)$  does \emph{not} lie in the multiplier ideal
$\Jj(\aa)$. But for all $c<1$, clearly $y \in \Jj(\aa^c)$.

To pave the way for clean proofs we need to formalize our setup slightly and
recall some results from toric geometry.

\subsection{Basic notions from toric geometry.}
Note that $k[X]=k[x_1,\ldots,x_n]$ carries a natural $\ZZ^n$-grading by giving
a monomial $x^m=x_1^{m_1}\cdot\ldots\cdot x_n^{m_n}$ degree $m \in \ZZ^n$.
Equivalently we note that the $n$-dimensional torus
\[ T^n \ =\ \Spec k[x_1^{\pm 1},\ldots,x_n^{\pm
1}] \ \cong \ (k^*)^n
\] acts on $k[X]$ via $\lambda \cdot x^m = \lambda^m
x^m$ for $\lambda \in (k^*)^n$. In terms of the varieties this means that
$X=\AA^n$ contains the torus $T^n$ as a dense open subset, and the action of
$T^n$ on itself naturally extends to an action of $T^n$ on all of $X$. Under
this action, the torus fixed ($=\ZZ^n$-graded) ideals are precisely the
monomial ideals. We denote the lattice $\ZZ^n$ in which the grading takes place
by $M$. It is just the lattice of the exponents of the Laurent monomials of
$k[T^n]$.

As indicated above, the Newton polytope $\Newt(\aa)$ of a monomial ideal $\aa$
is just the convex hull in $M_\RR=M \tensor_\ZZ \RR$ of the set $\{\, m \in M\,
|\, x^m \in \aa\,\}$. The Newton polytope of a principal ideal $(x^v)$ is just
the positive orthant in $M_\RR$ shifted by $v$. In general, the Newton polytope
of any ideal is an unbounded region contained in the first orthant. With every
point $v$ the Newton polytope also contains the first orthant shifted by $v$.

\begin{exercise}\label{ex.intclomon} Let $\aa$ be a
monomial ideal in $k[x_1,\ldots,x_n]$. Then the lattice points (viewed as
exponents) in the Newton polytope $\Newt(\aa)$ of $\aa$ define an ideal
$\bar{\aa} \supseteq \aa$. Show that $\bar{\aa}$ is the integral closure of
$\aa$ (see \cite{Fulton.Toric}).
\end{exercise}

The property of $X=\AA^n$ to contain the torus $T^n$ as a dense open set such
that the action of $T^n$ on itself extends to an action on $X$ as just
described is the definition of a \emph{toric variety}. The language of toric
varieties is the most natural to phrase, prove (and generalize, see
\cite{Bli.MultToric}) Howald's result. To set this up completely would take us
somewhat afield, so we choose to take a more direct approach using only a bare
minimum of toric geometry.

A first fact we have to take without proof from the theory of toric varieties
is that log resolutions of torus fixed ideals of $k[X]$ exist in the category
of toric varieties.\footnote{To be precise, a toric variety comes with the
datum of the torus embedding $T^n \subseteq X$. Maps of toric varieties are
such that they preserve the torus action.}
\begin{theorem}
    Let $\aa \subseteq k[x_1,\ldots,x_n]$ be a monomial ideal. Then there is a
    log resolution $\mu: Y \to X$ of $\aa$ such that $\mu$ is a map of toric
    varieties and consequently $\aa\cdot \Oo_Y = \Oo_Y(-F)$ is such that $F$ is fixed by the
    torus action on $Y$.
\end{theorem}
\begin{proof}[Indication of proof.]
This follows from the theory of toric varieties. First one takes the normalized
blowup of $\aa$, which is a (possibly singular) toric variety since $\aa$ was a
torus invariant ideal. Then one torically resolves the singularities of the
resulting variety as described in \cite{Fulton.Toric}. Note that this is a much
easier task than resolution of singularities in general. It comes down to a
purely combinatorial procedure.

An alternative proof could use Encinas and Villamayor's \cite{EncVill.Desing}
equivariant resolution of singularities. They give an algorithmic procedure of
constructing a log resolution of $\aa$ such that the torus action is preserved
--- that is by only blowing up along torus fixed centers.
\end{proof}

\subsubsection{Toric Divisors} A toric variety $X$ has a finite set of torus
fixed prime (Weil) divisors. Indeed, since an arbitrary torus fixed prime
divisor cannot meet the torus ($T^n$ acts transitive on itself and is dense in
$X$) it has to lie in the boundary $Y-T^n$, which is a variety of dimension
$\leq n-1$ and thus can only contain finitely many components of dimension
$n-1$. Furthermore, these torus fixed prime divisors $E_1,\ldots,E_r$ generate
the lattice of all torus fixed divisors which we shall denote by $L^X$. We
denote the sum of all torus invariant prime divisors $E_1+\ldots+E_r$ by $1_X$.

The torus invariant rational functions of a toric variety are just the Laurent
monomials $x_1^{m_1}\cdot\ldots\cdot x_n^{m_n} \in k[T^n]$. For the toric
variety $X=\AA^n$ one clearly has the identification of $M$, the lattice of
exponents, with $L^X$ by sending $m$ to $\div x^m$. In general this map will
not be surjective and its image is precisely the set of torus invariant Cartier
divisors. We note the following easy lemma which will nevertheless play an
important role in our proof of Theorem \ref{thm.mon}. It makes precise the idea
that a log resolution of a monomial ideal $\aa$ corresponds to turning its
Newton polytope $\Newt(\aa) \subseteq M_\RR$ into a translate of the first
orthant in $L^X_\RR$.
\begin{lemma}\label{lem.techmon}
    Let $\mu: Y \to X=\Spec k[x_1,\ldots,x_n]$ be a toric resolution of the
    monomial ideal $\aa \subseteq k[x_1,\ldots,x_n]$ such that $\aa \cdot \Oo_Y
    = \Oo_Y(-F)$. Then, for $m \in M$ we have
    \[
        c \cdot m \in c'\Newt(\aa) \iff c\cdot\mu^* \div x^m \geq c' \cdot F
    \]
    for all rational $c,c' > 0$.
\end{lemma}
\begin{proof}
    We first show the case $c=c'=1$. Assume that $m \in \Newt(\aa)$. By
    Exercise \ref{ex.intclomon}, this is equivalent to $x^m \in \bar{\aa}$, the
    integral closure of $\aa$. Since, by Exercise \ref{ex.intclo}, $\bar{\aa} =
    \mu_*\Oo_Y(-F)$ it follows that $x^m \in \bar{\aa}$
    if and only if $\mu^*x^m \in \Oo_Y(-F)$. This, finally, is equivalent
    to $\mu^*(\div x^m) \geq F$.

    For the general case express $c$ and $c'$ as integer fractions.
    Clearing denominators and noticing that for an integer $a$ one has $a
    \Newt(\aa)=\Newt(\aa^a)$ one reduces to the previous case.
\end{proof}
\subsubsection{Canonical divisor} As the final ingredient for computing the
multiplier ideal we need an understanding of the canonical divisor (class) of a
toric variety.
\begin{lemma}\label{lem.can}
    Let $X$ be a (smooth) toric variety and let $E_1,\ldots,E_r$ denote the collection
    of all torus invariant prime Weil divisors. Then the canonical divisor is
    $K_Y=-\sum E_i=-1_X$.
\end{lemma}
We leave the proof as an exercise or alternatively refer to \cite{Fulton.Toric}
or \cite{Danilov78} for this basic result. We verify it for $X=\AA^n$. Then
$E_i=(x_i=0)$ for $i=1,\ldots,n$ are the torus invariant divisors and $K_X$ is
represented by the divisor of the $T^n$-invariant rational $n$-form
$\tfrac{dx_1}{x_1}\wedge\ldots\wedge\tfrac{dx_n}{x_n}$, which is
$-(E_1+\ldots+E_n)$. As a consequence of the last lemma we get the following
lemma.
\begin{lemma}\label{lem.mon1X}
    Let $\mu: Y \to X=\AA^n$ be a birational map of (smooth) toric varieties. Then
    $K_{Y/X}=\mu^* 1_X - 1_Y$ and the support of $\mu^* 1_X$ is equal to the
    support of $1_Y$.
\end{lemma}
\begin{proof}
    As the strict transform of a torus invariant divisor on $X$ is
    a torus invariant divisor on $Y$ it follows that $\mu^* 1_X - 1_Y$ is
    supported on the exceptional locus of $\mu$. Since $-1_X$ represents the
    canonical class $K_X$ and respectively for $Y$, the first assertion follows
    from the definition of $K_{Y/X}$. Since $\mu^* 1_X$ is torus invariant clearly
    its support is included in $1_Y$. Since $\mu$ is an isomorphism over the
    torus $T^n \subseteq X$ it follows that $\mu^{-1}(1_X) \supseteq 1_Y$ which implies the second
    assertion.
\end{proof}

\begin{exercise}
    This exercise shows how to avoid taking Lemma \ref{lem.can} on faith but
    instead using a result of Russel Goward \cite{Goward.PrincMon} which states that a log resolution
    of a monomial ideal can be obtained by a sequence of monomial blowups.

    A \emph{monomial blowup} $Y=\Bl_Z(Y)$ of $\AA^n$ is the blowing up
    of $\AA^n$ at the intersection $Z$ of some of the coordinate hyperplanes
    $E_i=(x_i=0)$ of $\AA^n$.

    For such monomial blowup $\mu: Y=\Bl_Z(X) \to X \cong \AA^n$ show that $Y$
    is a smooth toric variety which is canonically covered by $\codim(Z,X)$ many
    $\AA^n$ patches. Show that $1_Y=E_1+\ldots+E_n+E$ where $E$ is the exceptional
    divisor of $\mu$. Via a direct calculation verify the assertions of the
    last two lemmata for $Y$.

    Since a monomial blowup is canonically covered by affine spaces one can
    repeat the process and obtains the notion of a \emph{sequence of monomial
    blowups}. Using Goward's result show directly that a monomial ideal has a toric log
    resolution $\mu: Y \to \AA$ with the properties as in Lemma \ref{lem.mon1X}.
\end{exercise}

\subsection{Proof of Theorem \ref{thm.mon}}
By the existence of a toric (or equivariant) log resolution of a monomial ideal
$\aa$ it follows immediately that the multiplier ideal $\Jj(\aa^c)$ is also
generated by monomials. Thus, in order to determine $\Jj(\aa^c)$ it is enough
to decide which monomials $x^m$ lie in $\Jj(\aa^c)$. With our preparations this
now an easy task.
\begin{proof}[Proof of Theorem \ref{thm.mon}]
    As usual we denote $\Spec k[x_1,\ldots,x_n]$ by $X$ and let $\mu:Y \to X$
    be a toric log resolution of $\aa$ such that $\aa \cdot \Oo_Y =
    \Oo_Y(-F)$.

    Abusing notation by identifying $\div(x_1\cdot\ldots\cdot x_n) = 1_X \in L^X$ with $(1,\ldots,1) \in M$,
    the condition of the theorem that $m + 1_X$ is in the interior of the Newton polytope
    $c\cdot\Newt(\aa)$ is equivalent to
    \[
        m + 1_X - \epsilon 1_X \in c\Newt(\aa)
    \]
    for small enough rational $\epsilon > 0$. By Lemma
    \ref{lem.techmon} this holds if and only if
    \[
        \mu^* \div g + \mu^* 1_X - \epsilon \mu^*1_X
\geq c\,F.
    \]
    Using the formula $K_{Y/X}=\mu^*1_X-  1_Y$ from
Lemma
    \ref{lem.can} this is furthermore equivalent to
    \[
        \mu^*\div g + K_{Y/X}+\floor{1_{Y} -
\epsilon\mu^*1_X -c\,F} \geq 0
    \]
    for sufficiently small $\epsilon > 0$.
    Since by Lemma \ref{lem.mon1X}, $\mu^*1_X$ is effective with the same support
    as $1_{Y}$ it follows that all coefficients appearing
    in $1_{Y} -\epsilon\mu^*1_X$ are very close to but
    strictly smaller than 1 for
    small $\epsilon > 0$. Therefore, $\floor{1_{Y} -
\epsilon\mu^*1_X -cF}
    = \ceil{-cF} = -\floor{cF}$.
    Thus we can finish our chain of equivalences with
    \[
        \mu^* \div g \geq -K_{Y/X} + \floor{cF}
    \]
    which says nothing but that $g \in \Jj(\aa^c)$.
\end{proof}

This formula for the multiplier ideal of a monomial ideal will be applied in
the next section to concretely compute certain invariants arising from
multiplier ideals.

\section{Invariants arising from multiplier ideals and
applications}\label{sec.invariants}

We keep the notation of a smooth  affine variety $X$
over an algebraically closed field of characteristic
zero, and an ideal
$\aa
\subseteq
k[X]$. In this section we use multiplier ideals to
attach some invariants to $\fra$, and we study their
influence on some algebraic questions.

\subsection{The log canonical threshold} If $c > 0$ is
very small, then $\MI{\fra^c} = k[X]$.  For large
$c$, on the other hand, the multiplier ideal
$\Jj(\aa^c)$ is clearly nontrivial. This leads one to
define:
\begin{definition} The \emph{log canonical threshold}
of $\fra$ is the number
\[
    \lct (\aa) = \lct (X,\aa) = \inf\{\, c > 0\, |\,
\Jj(\aa^c)
    \neq
    \Oo_X \, \}.
\]
\end{definition}
The following exercise shows that $\lct(\fra)$ is a rational number, and that
the infimum appearing in the definition is actually a minimum. Consequently,
the log canonical threshold is just the smallest $c > 0$ such that $\Jj(\aa^c)$
is nontrivial.
\begin{exercise} \label{LCT.via.discrepencies}
    As usual, fixing notation of a log resolution $\mu:
Y \to X$ with
    $\aa\cdot\Oo_{Y}=\sum r_iE_i$ and $K_{Y/X}=\sum
    b_iE_i$, show that $\lct(X,\aa)= \min \{\,
    \frac{b_i+1}{r_i}\,\}$.
\end{exercise}   Recall the notions from singularity theory
\cite{Kollar.Sings.of.Pairs} in which a pair $(X,\aa^{c})$ is called
\emph{log terminal} if and only if $b_i - cr_i + 1
> 0$ for all $i$. It is called \emph{log canonical} if and only if $b_i -cr_i +
1\geq 0$ for all $i$. The last exercise also shows that $(X,\aa^c)$ is log
terminal if and only if the multiplier ideal $\Jj(\aa^c)$ is trivial.

\begin{example} Continuing previous examples we observe
that $\lct((x^2,y^2))=1$ and $\lct((x^2,y^3))=\frac{5}{6}$.
\end{example}

\begin{example}[The log canonical threshold of a
monomial ideal] The formula for the multiplier ideal of a monomial ideal $\aa$
on $X = \Spec k[x_1,\ldots,x_n]$ shows that $\Jj(a^c)$ is trivial if and only
if $1_X=(1,\ldots,1)$ is in the interior of the Newton polytope
$c\,\Newt(\aa)$. This allows to compute the log canonical threshold of $\aa$:
$\lct(\aa)$ is the largest $t>0$ such that $1_X \in t \cdot \Newt(\aa)$.
\end{example}

\begin{example} As a special case of the previous
example, take
\[ \aa\ =\ (x_1^{a_1},\ldots,x_n^{a_n}). \] Then
the Newton polytope is the subset of the first orthant consisting of points
$(v_1,\ldots,v_n)$ satisfying $\sum \frac{v_i}{a_i} \geq 1$. Therefore $1_X \in
t\cdot \Newt(\aa)$ if and only if $\sum \frac{1}{a_i} \geq t$. In particular,
$\lct(\aa) = \sum \frac{1}{a_i}$.
\end{example}

\subsection{Jumping numbers} The
log-canonical threshold measures the triviality or non-triviality of a
multiplier ideal. By using the full algebraic structure of these ideals, it is
natural to see this threshold as merely the first of a sequence of invariants.
These so-called jumping numbers were first considered (at least implicitly) in
work of Libgober, Loeser and Vaqui\'e (\cite{Libgober},
\cite{Loeser.Vaquie}). They are studied more systematically in the paper
\cite{ELSV}.

We start with a lemma:
\begin{lemma}\label{prop.jump}
    For $\aa \subseteq \Oo_X$, there is an increasing discrete
sequence of
    rational numbers
    \[
        0 \ = \ \xi_0\ <\  \xi_1\  < \ \xi_2 \ <\
\ldots
    \]
    such that $\Jj(\aa^c)$ is constant for $\xi_i\leq c
<
    \xi_{i+1}$ and $\Jj(\aa^{\xi_i})
\varsupsetneq\Jj(\aa^{\xi_{i+1}})$.
\end{lemma}
\noi We leave the (easy) proof to the reader.

The $\xi_i = \xi_i(\fra)$
are called the
\emph{jumping numbers} or \textit{jumping
coefficients} of
$\aa$. Referring to the log resolution $\mu$ appearing
in Example \ref{LCT.via.discrepencies}, note that the
only candidates for jumping numbers are those
$c$ such that
$cr_i$ is an integer for some $i$. Clearly the first
jumping number
$\xi_1(\aa)$ is the log canonical threshold $\lct(\aa)$.

\begin{example}[Jumping numbers of monomial ideals]
Let $\fra \subseteq k[x_1, \ldots, x_n]$ be a monomial ideal. For the
multiplier ideal $\Jj(\aa^c)$ to jump at $c=\xi$, it is equivalent that some
monomial, say $x^v$, is in $\Jj(\aa^\xi)$ but not in $\Jj(\aa^{\xi-\epsilon})$
for all $\epsilon > 0$. Thus, the largest $\xi>0$ such that $v+(1,\ldots,1) \in
\xi \Newt(\aa)$ is a jumping number. Doing this construction for all $v \in
\NN^n$ one obtains all jumping numbers of $\aa$ (this uses the fact that the
multiplier ideal of a monomial ideal is a monomial ideal).
\end{example}

\begin{exercise} Consider again $\aa =
(x_1^{a_1},\ldots,x_n^{a_n})$. Then the jumping numbers
of $\fra$  are precisely the
rational numbers of the form
\[
    \tfrac{v_1+1}{a_1}+\ldots+\tfrac{v_n+1}{a_n}
\] where $(v_1,\ldots,v_n)$ ranges over $\NN^n$. Note
however that different vectors $(v_1, \ldots,v_n)$ may
give the same jumping number.
\end{exercise}

\begin{figure}
\includegraphics[scale = .8]{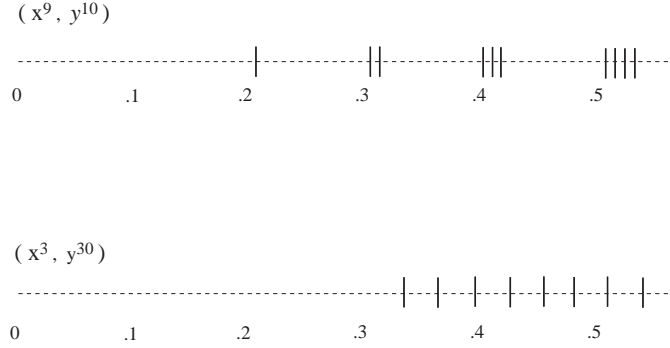}
\caption{Jumping numbers of $(x^9, y^{10})$ and $(x^3,
y^{30})$}\label{Jumping.Nos.Picture}
\end{figure}

It is instructive to picture the jumping numbers of an
ideal graphically.
 Figure \ref{Jumping.Nos.Picture},
taken from
\cite{ELSV}, shows the jumping numbers of the two
ideals $(x^9, y^{10})$ and $(x^3, y^{30})$: the
exponents are chosen so that the two ideals have the
same Samuel multiplicity, and so that the pictured
jumping coefficients occur ``with multiplicity one" (in
a sense whose meaning we leave to the reader).

\subsection{Jumping length}
Jumping numbers give rise to an additional invariant in
the case of principal ideals.
\begin{lemma} \label{Jump.Nos.Poly.Lemma}
    Let $f\in k[X]$ be a non-zero function. Then $\Jj(f)
= (f)$ but for
$c<1$ one has
    $(f) \varsubsetneq \Jj(f^c)$. In other words,  $\xi =
1$ is a
    jumping number of the principal ideal $(f)$.
\end{lemma}
Deferring the proof for a moment, we note that the Lemma means that
$\xi_{\ell}(f) = 1$ for some index $\ell$. We define $\ell = \ell(f)$ to be the
\textit{jumping length} of $f$. Thus $\ell(f)$ counts the number of jumping
coefficients of $(f)$ that are $\le 1$.

\begin{example}
    Let $f = x^4 + y^3 \in \CC[x,y]$. One can show that
$f$ is sufficiently
    generic so that  $\Jj(f^c) =
    \Jj((x^4,y^3)^c)$ provided that $c < 1$. Therefore
the first few jumping numbers of $f$ are
       \[
        0\ < \ \lct (f)\  = \ \tfrac{1}{4} +
\tfrac{1}{3}\  <\
\tfrac{2}{4} \ +\
        \tfrac{1}{3}\  <\  \tfrac{1}{4} + \tfrac{2}{3}
\ <\  1,
    \]
    and $\ell(f) = 4$.
\end{example}

\begin{proof} [Proof of Lemma \ref{Jump.Nos.Poly.Lemma}]
    Let $\mu: Y \to X$ be a log resolution of $(f)$ and
denote the integral divisor
    $(f=0)$ by $D = \sum a_i D_i$. Clearly, $\aa \cdot
\Oo_Y =
\Oo_Y(-\mu^*D)$ and
$\mu^*D$ is
    also an integral divisor. Thus
    \[
    \begin{split}
        \Jj(f)&=\mu_*\Oo_Y(K_{Y/X}-\mu^*D) \\ &=
        \mu_*(\Oo_Y(K_{Y/X})\tensor\mu^*\Oo_X(-D)) \\
&=\Oo_X \tensor
\Oo_X(-D)\\ &=(f).
    \end{split}
    \]
On the other hand, choose a general point $x \in D_i$
on any of the components of $D = \text{div}(f) = \sum
a_i D_i$. Then
$\mu$ is an isomorphism over $x$ and consequently
\[ \text{ord}_{D_i} \big ( \MI{ f^c} \big ) \ < \ a_i \
\text{ for }  0 < c < 1. \]
Therefore $\MI{(f)^c} \subsetneqq (f)$ whenever $c <
1$.
\end{proof}

Finally, we note that the jumping length can be related to other invariants of
the singularities of $f$:
\begin{proposition}[\cite{ELSV}] \label{Tyurina.No.Bound}
    Assume the hypersurface defined by the vanishing of $f$ has at worst an isolated singularity
at $x \in X$.
    Then
    \[
        \ell(f) \leq \tau(f,x) + 1,
    \]
    where $\tau(f,x)$ is the Tjurina number of
$f$ at $x$,  defined as the colength in $\Oo_{x,X}$ of $(f,\frac{\partial
    f}{\partial z_1},\ldots,\frac{\partial
    f}{\partial z_n})$ for $z_1,\ldots,z_n$ parameters
around $x$.
\end{proposition}

\subsection{Application to uniform Artin-Rees numbers}
We next discuss  a result relating jumping lengths to
uniform Artin-Rees numbers of a principal ideal.

To set the stage, recall the statement of the Artin-Rees lemma in a simple
setting:
\begin{theorem*}[Artin-Rees]
    Let $\bb$ be an ideal and $f$ an element of $k[X]$. There exists
    an integer $k = k(f,\bb)$ such that
    \[
        \bb^m \cap (f) \ \subseteq \ \bb^{m-k}\cdot(f)
    \]
    for all $m \geq k$. In other words, if $fg \in
\bb^m$ then  $g
    \in \bb^{m-k}$.
\end{theorem*}
\noi Classically, $k$ is allowed to depend both on $\frb$ and $f$. However in
his paper \cite{Huneke.UniformBounds}, Huneke showed that in fact there is a
single integer $k=k(f)$ which works simultaneously for all ideals
$\bb$.\footnote{We stress that both the classical Artin-Rees Lemma and Huneke's
theorem are valid in a much more general setting.} Any such $k$ is called a
\textit{uniform Artin-Rees number} of $f$.

The next result shows that the jumping length gives an
effective estimate (of moderate size!) for uniform
Artin-Rees numbers.
\begin{theorem}[\cite{ELSV}]\label{thm.ArtReesNum}
    As above, write $\ell(f)$ for the jumping length
of $f$. Then the integer
$k=\ell(f)\cdot
\dim X$    is a uniform Artin-Rees number of $f$.
\end{theorem}
If $f$ defines a smooth hypersurface, its jumping length is $1$ and it follows
that $n=\dim X$ is a uniform Artin-Rees number in this case. (In fact, Huneke
showed that $n - 1$ also works in this case.)

If $f$ defines a hypersurface with only an isolated singular point $x\in X$, it
follows from Proposition \ref{Tyurina.No.Bound} and the Theorem that $k = n
\cdot \big(\tau(f,x) + 1) \big)$ is a uniform Artin-Rees number. (One can show
using the next Lemma and some observations of Huneke that in fact $k =
\tau(f,x) + n$ also works: see \cite[\S 3]{ELSV}.)

The essential input to Theorem \ref{thm.ArtReesNum} is
a statement involving  consecutive jumping coefficients:
\begin{lemma} \label{Lemma.UAR.Pr.Ideals}
Consider  two consecutive jumping numbers
\[\xi \, =\,
\xi_i(f) \ <\
\xi_{i+1}(f) \, =\, \xi'\]  of $f$, and let  $\bb
\subseteq k[X]$ be any ideal. Then given a natural
number
$m>n=\dim X$, one has
    \[
        \bb^m\cdot\Jj(f^{\xi})\, \cap \, \Jj(f^{\xi'})
\ \subseteq\
        \bb^{m-n}\cdot\Jj(f^{\xi'}).
    \]
\end{lemma}
\noi We will deduce this from Skoda's theorem in the
next section. In the meantime, we observe that it leads
immediately to the
\begin{proof} [Proof of Theorem \ref{thm.ArtReesNum}]
We apply the previous Lemma  repetitively to successive
jumping numbers in the chain of multiplier ideals

\[
    k[X] \ = \ \Jj(f^0) \ \varsupsetneq \Jj(f^{\xi_1})
\ \varsupsetneq\
    \Jj(f^{\xi_2})\ \varsupsetneq \  \ldots \ \varsupsetneq
\
\Jj(f^{\xi_\ell}) =
    \Jj(f) = (f).
\] After  further intersection with $(f)$ one finds:
\begin{align*}
    \bb^m   \cap (f) \ &\subseteq \
\bb^{m-n} \cdot \Jj(f^{\xi_1}) \cap (f) \\
&\subseteq
    \bb^{m-2n}\cdot \Jj(f^{\xi_2}) \cap (f) \\
& \quad
\ldots \\
&\subseteq
    \bb^{m-\ell n}(f),
\end{align*}as required.
\end{proof}

\begin{remark} When $\fra = (f )$ is a principal ideal, the jumping
numbers of $f$ are related to other invariants appearing in the
literature. In particular, if $f$ has an isolated singularity then
(suitable translates of) the jumping coefficients appear in the
Hodge-theoretically defined \textit{spectrum} of $f$. See \cite[\S
5]{ELSV} for precise statements and references.
\end{remark}

\section{Further Local Properties Of Multiplier Ideals}
In this section we discuss some results involving the local behavior of
multiplier ideals. We start with  Skoda's theorem and some variants. Then we
discuss the restriction and subadditivity theorems, which will be used in the
next section.

\subsection{Skoda's theorem} An important (and early)
example of a uniform result in local algebra was established by Skoda and
Brian\c con \cite{BS.74} using analytic results of Skoda \cite{Skoda72}.
In our language, Skoda's result is this:

\begin{theorem}[Skoda's Theorem, I]\label{thm.version1}
    Consider any ideal $\bb \subseteq k[X]$ with $X$
smooth of dimension
    $n$. Then for all $m \geq n$
    \[
        \Jj(\bb^m) \ = \ \bb \cdot \Jj(\bb^{m-1})\  =\
\ldots\ =\
        \bb^{m+1-n}\cdot \Jj(\bb^{n-1}).
    \]
\end{theorem}

\begin{remark} As Hochster noted in his lectures, the
statement in \cite{Skoda72} has a more analytic flavor.
In fact, using the analytic interpretation of
multiplier ideals (\S \ref{Analytic.Interp.Subsectn})
one sees that (the analytic analogue of)  Theorem
\ref{thm.version1} is essentially equivalent to the
following statement.
\begin{quote}
Suppose that
 $\bb$ is generated by
$(g_1,\ldots,g_t)$, and that $f$ is a holomorphic
function such that
\begin{equation}
 \int
\frac{\abs{f}^2}{(\sum
     \abs{g_i}^{2})^{m}} \ < \ \infty   \notag
\end{equation} for some $m \ge n = \dim X$.   Then
locally there exist holomorphic functions $h_i$ such
that
$ f=\sum h_ig_i $, and moreover each of the $h_i$
satisfies the local integrability condition $\int
\frac{\vert h_i
\vert^2}{(\sum  \vert g_i \vert^2)^{m-1}} <
\infty$.\end{quote}  (The hypothesis expresses
the membership of $f$ in
$\MI{\frb^m}^{\text{an}}$ and the conclusion writes $f$
as belonging to $\frb^{\text{an}}
\cdot
\MI{\frb^{m-1}}^{\text{an}}$.)
\end{remark}

As a corollary of Skoda's theorem, one obtains the
classical theorem of Brian\c{c}on-Skoda.
\begin{corollary}[Brian\c{c}on-Skoda]
    With the notation as before,
    \[
        \bar{\bb^m}\ \subseteq\ \Jj(\bb^m) \ \subseteq\
\bb^{m+1-n}
    \]
    where $\bar{\phantom{\bb}}$ denotes the integral
closure and $n =
\dim X$.
\end{corollary}

\begin{proof}[Sketch of proof of Theorem \ref{thm.version1}] The argument
follows ideas of Teissier and  Lipman. We choose
generators
$g_1,\ldots,g_k$ for the ideal $\bb$ and fix a log resolution $\mu: Y \to
X$ of $\bb$ with $\bb
\cdot \Oo_{Y} = \Oo_{Y}(-F)$. Write $g^{\prime}_i=\mu^*(g_i) \in
\Gamma(Y,\Oo_{Y}(-F))$ to define the surjective map
\begin{equation}\label{eqn.FirstMapKoszul}
    \oplus_{i=0}^{k} \Oo_{Y} \to \Oo_{Y}(-F)
\end{equation}
by sending $(x_1,\ldots,x_k)$ to $\sum x_ig^{\prime}_i$. Tensoring this map
with $\Oo_{Y}(K_{Y/X}-(m-1)F)$ yields the surjection
\[
    \oplus_{i=1}^k \Oo_{Y}(K_{Y/X}-(m-1)F) \to[\phi]
    \Oo_{Y}(K_{Y/X}-mF).
\]
Further applying $\mu_*$ we get the map $\oplus_{i=0}^k \Jj(\bb^{m-1})
\to[\mu_*\phi] \Jj(\bb^m)$ which again sends a tuple $(y_1, \ldots, y_k)$ to
$\sum y_ig_i$. Therefore, the image of $\mu_*(\phi)$ is
\[
    \op{Image}(\mu_*\phi) \ = \ \bb\Jj(\bb^{m-1})
\ \subseteq\ \Jj(\bb^m).
\]
What remains to show is that $\mu_*\phi$ is surjective. For this consider the
Koszul complex  on the $g_i'$ on $Y$ which resolves the map in
(\ref{eqn.FirstMapKoszul}).
\[
\begin{split}
    0 \to \Oo_Y((k-1)F) &\to \oplus^{k} \Oo_Y((k-2)F) \to \ldots
    \\ \ldots & \to \oplus^{\binom{k}{2}} \Oo_Y(F) \to \oplus^k \Oo_Y \to
    \Oo_Y(-F) \to 0.
\end{split}
\]
As above, tensor through by $\Oo_Y(K_{Y/X}-(m-1)F)$ to get a resolution of
$\phi$. Local vanishing (Theorem \ref{Local.Vanishing.Thm}) applies to the $m
\geq n = \dim X$ terms on the right. Chasing through the sequence while taking
direct images then gives the required surjectivity. See \cite[Chapter
9]{PAG} or \cite{ELNull} for details.
\end{proof}

 It will be useful to have a variant involving several
ideals and fractional coefficients.  For this we extend
slightly  the definition of multiplier ideals.

\subsubsection{Mixed multiplier ideals} Fix a sequence
of ideals
$\aa_1,\ldots,\aa_t$ and positive rational numbers $c_1,
\ldots, c_t$. Then we define the multiplier ideal
\[
    \Jj(\aa_1^{c_1}\cdot\ldots\cdot\aa_t^{c_t})
\] starting with a log
resolution $\mu: Y
\to X$ of the product $\aa_1 \cdot \ldots \cdot \aa_t$.
Since this is at the same time also a log resolution of
each $\aa_i$ write $\aa_i\cdot
\Oo_{Y} =
\Oo_{Y}(-F_i)$ for simple normal crossing divisors
$F_i$.
\begin{definition}
    With the notation as indicated, the mixed
multiplier ideal is
    \[
\Jj(\aa_1^{c_1}\cdot\ldots\cdot\aa_t^{c_t} ) \ =
\
\mu_*(\Oo_{Y}(K_{Y/X}-\floor{c_1F_1+\ldots+c_tF_t})).
    \]
    As before, this definition is independent of the
chosen log resolution.
\end{definition} Note that once again we do not
attempt to assign any meaning to the expression
$\aa_1^{c_1}\cdot\ldots\cdot\aa_t^{c_t}$ in the
argument of $\Jj$. This expression is meaningful a
priori whenever all $c_i$ are positive integers and our
definition is consistent with this prior meaning.

With this generalization of the concept of multiplier ideals we get the
following variant of Skoda's theorem.
\begin{theorem}[Skoda's Theorem, II]
\label{Skoda.Thm.II}
    For every integer $c \ge n = \dim X$ and any $d > 0$ one has
    \[
        \Jj(\aa_1^{c}\cdot\aa_2^{d}) \ =\
        \aa_1^{c-(n-1)}\Jj(\aa_1^{n-1}\cdot\aa_2^{d}).
    \]
\end{theorem} The proof of this result is only a
technical complication of the proof of the First Version we discussed above. We
refer to \cite[Chapter 9]{PAG} for details.

We conclude by using Skoda's Theorem to prove (a slight
generalization of) the Lemma \ref{Lemma.UAR.Pr.Ideals}
underlying the results on uniform Artin-Rees numbers
in the previous section.
\begin{lemma}\label{lem.ArtResUsed}
    Let $\aa \subseteq k[X]$ be an ideal and let $\xi
< \xi'$ be
    consecutive jumping numbers of $\aa$. Then for $m >
n$ we have
    \[
        \bb^m\cdot \Jj(\aa^{\xi})\ \cap \
\Jj(\aa^{\xi'})
 \ \subseteq \
        \bb^{m-n}\cdot \Jj(\aa^{\xi'})
    \]
    for all ideals $\bb \subseteq k[X]$.
\end{lemma}
\begin{proof}
    We first claim that
    \[
        \bb^m\Jj(\aa^{\xi}) \, \cap \, \Jj(\aa^{\xi'})
       \  \subseteq \ \Jj(\bb^{m-1} \cdot \aa^{\xi'}).
    \]
This is shown via a simple computation.  In fact, to begin with  one can
replace $\xi$ by $c \in [\xi, \xi')$ arbitrarily close to $\xi'$ since this
does not change the statement. Let $\mu:Y \to X$ be a common log resolution of
$\aa$ and $\bb$ such that $\aa\cdot\Oo_Y = \Oo_Y(-A)$ and $\bb\cdot\Oo_Y =
\Oo_Y(-B)$. Let $E$ be a prime divisor on $Y$ and denote by $a$, $b$ and $e$
the coefficient of $E$ in $A$, $B$ and $K_{Y/X}$, respectively. Then $f$ is in
the left-hand side if and only if
    \[
        \ord_E f \ \geq\ \max(-e+mb+\floor{ca}\, ,\, -e+\floor{\xi'a}).
    \]
    If $b=0$ this implies that $\ord_E f \geq -e+(m-1)b+\floor{\xi'a}$. If
    $b\neq 0$ then $b$ is a positive integer $\geq 1$. Since $c$ is
arbitrarily  close
    to $\xi'$ we get $\floor{\xi' a} - b \leq \floor{\xi' a} - 1 \leq \floor{ca}$.
    Adding $-e+mb$ it follows that also in this case $\ord_E f \geq
    -e+(m-1)b+\floor{\xi'a}$. Since this holds for all $E$ it follows that
    $f\in\Jj(\bb^{m-1} \cdot \aa^{\xi'})$.

    Now, using Theorem \ref{Skoda.Thm.II} we deduce
   \[
        \Jj(\bb^{m-1}\cdot\aa^{\xi'}) \ \subseteq\
        \bb^{m-n}\Jj(\bb^{n-1}\cdot\aa^{\xi'})
\ \subseteq \
        \bb^{m-n}\Jj(\aa^{\xi'}).
    \]
  Putting all the inclusions together, the Lemma
follows.
\end{proof}
\begin{exercise} \label{Periodicityu.Jump.Nos}
    Let $\aa \subseteq k[X]$ be an ideal. Starting at
$\dim X-1$, the
    jumping numbers are periodic with period 1. That
is, $\xi
    \geq \dim X-1$ is a jumping number if and only if
$\xi+1$ is a jumping number.
\end{exercise}

\subsection{Restriction theorem}

The next result deals with restrictions of multiplier ideals. Consider a smooth
subvariety $Y \subseteq X$ and an ideal $\frb \subseteq k[X]$ which does not
vanish on $Y$. There are then two ways to get an ideal on $Y$. First, one can
compute the multiplier ideal $\MI{X, \frb^c}$ on $X$ and then restrict it
to
$Y$. On the other hand, one can also restrict $\frb$ to $Y$ and then compute
the multiplier ideal on $Y$ of this restricted ideal.   The Restriction
Theorem -- which is arguably the most important local property of
multiplier ideals -- states that there is always an inclusion among these
ideals on $Y$.
\begin{theorem}[Restriction Theorem]\label{thm.restr}
    Let $Y \subseteq X$ be a smooth subvariety of $X$ and
$\bb$ an ideal of
    $k[X]$ such that $Y$ is not contained in the zero
locus of $\bb$. Then
    \[
        \Jj(Y,(\bb\cdot k[Y])^c)\  \subseteq\
\Jj(X,\bb^c)\cdot k[Y].
    \]
\end{theorem}
\noi One can think of the Theorem as reflecting the
principle that singularities can only get worse under
restriction.

In the present setting, the result is due to Esnault and Viehweg
\cite[Proposition 7.5]{EV}
 When $Y$ is a hypersurface, the statement is proved
using the local vanishing theorem \ref{Local.Vanishing.Thm}. Since in any event
a smooth subvariety is a local complete intersection, the general case then
follows from this.

\begin{exercise}
Give an example where strict inclusion holds in the
Theorem.
\end{exercise}

\subsection{Subadditivity theorem} We conclude with a
result due to Demailly, Ein and the second author \cite{DEL}
concerning the multiplicative behavior of multiplier
ideals. This subadditivity theorem will be used in the
next section to obtain some uniform bounds on symbolic
powers of ideals.
\begin{theorem}[Subadditivity]\label{thm.subadd}
    Let $\aa$ and $\bb$ be ideals in $k[X]$. Then for
all $c,d > 0$ one has
    \[
        \Jj(\aa^c\cdot \bb^d) \ \subseteq\
\Jj(\aa^c)\cdot \Jj(\bb^d).
    \]
    In particular, for every positive integer $m$,
$\Jj(\aa^{cm})\subseteq
    \Jj(\aa^c)^m$.
\end{theorem}
\begin{proof}[Sketch of proof]\newcommand{\pr}{\op{p}}
    The idea of the proof is to pull back the data to
    the product $X \times X$ and then to restrict to the diagonal $\Delta$.
    Specifically, assume for simplicity that $c=d=1$, and consider the
    product
    \[
        \xymatrix{
        & {X\times X} \ar^{\pr_1}[dl]\ar_{\pr_2}[dr]&\\
        X &&X }
    \]
    along with its projections as indicated. For log
resolutions $\mu_1$ and
    $\mu_2$ of $\aa$ and $\bb$ respectively one can
verify that $\mu_1
\times
    \mu_2$ is a log resolution of the ideal
$\pr_1^{-1}(\aa) \cdot
\pr_2^{-1}(\bb)$
    on $X \times X$. Using this one shows that
    \[
        \Jj(X\times X,\pr_1^{-1}(\aa) \cdot
\pr_2^{-1}(\bb))\ = \
        \pr_1^{-1}\Jj(X,\aa)\cdot\pr_2^{-1}\Jj(X,\bb).
    \]
    Now let $\Delta \subseteq X\times X$ be the
diagonal. Apply the Restriction
    Theorem \ref{thm.restr} with $Y=\Delta$ to conclude
    \[
    \begin{split}
        \Jj(X,\aa\cdot\bb) \
&= \ \Jj(\Delta,\pr_1^{-1}(\aa)
\cdot
        \pr_2^{-1}(\bb)\cdot \Oo_{\Delta}) \\
            &\subseteq \ \Jj(X \times X, \pr_1^{-1}(\aa)
\cdot
            \pr_2^{-1}(\bb))\cdot\Oo_{\Delta} \\
&= \ \Jj(X,\aa) \cdot \Jj(X,\bb),
    \end{split}
    \]
   as required.\end{proof}

\section{Asymptotic Constructions} There are many
natural situations in geometry and algebra where one is forced to confront
rings or algebras that fail to be finitely generated. For example, if $D$
is a non-ample divisor on a projective variety $V$, then the section ring
$R(V,D) =
\oplus \Gamma(V,
\Oo_V(mD))$ is typically not finitely generated. Or likewise, if $\frq $ is a
radical ideal in some ring, the symbolic blow-up algebra $\oplus \frq^{(m)}$
likewise fails to be finitely generated in general. It is nonetheless possible
to extend the theory of multiplier ideals to such settings. It turns out that
there is
 finiteness built into the resulting multiplier
ideals that may not be present in the underlying
geometry or algebra. This has led to some of the most
interesting applications of the theory.

In the geometric setting, the asymptotic constructions
have been known for some time, but it was only with
Siu's work \cite{Siu98a} on deformation-invariance of
plurigenera that their power became clear. Here we
focus on an algebraic formulation of  the theory from
\cite{ELS1}. As before, we work with a smooth affine
variety $X$ defined over an algebraically closed field
$k$ of characteristic zero.

\subsection{Graded systems of ideals}
We start by defining certain collections of ideals, to
which we will later attach multiplier ideals.
\begin{definition}
    A \emph{graded system} or \textit{graded family of
ideals} is a family
$\fra_{\bullet} = \set{\aa_k}_{k \in
\NN}$ of
    ideals in $k[X]$ such that \[
\fra_{\ell} \cdot \aa_m \ \subseteq
\aa_{\ell+m}\] for all
    $\ell,m \geq 1$. To avoid trivialities, we also
assume that $\fra_k \ne (0) $ for $k \gg 1$.
\end{definition}
\noi The condition in the definition means that the direct sum
\[ R(\fra_{\bullet}) \ \overset{\text{def}}{=} \ k[X] \,
\oplus
  \fra_1
\oplus   \fra_2   \oplus  \ldots \] naturally carries a graded $k[X]$-algebra
structure and $R(\aa_\bullet)$ is called the \emph{Rees algebra} of
$\aa_\bullet$. In the interesting situations $R(\fra_{\bullet})$ is not
finitely generated, and it is here that the constructions of the present
section give something new. One can view graded systems as local objects
displaying complexities similar to those that arise from linear series on a
projective variety $V$.\footnote{If $D$ is an effective divisor on $V$, the
base ideals $\frb_k = \frb(|kD|) \subseteq \Oo_V$ form a graded family of ideal
sheaves on $V$: this is the prototypical example.}

\begin{example} We give several examples of graded
systems.
  \begin{enumerate} \label{First.Ex.GSI}
\item[(i).]   Let $\bb \subseteq k[X]$ be a fixed
ideal,  and set
$\aa_k=\bb^k$. One should view the resulting  as a
trivial example.
\item[(ii).]
 Let $Z \subseteq X$ be a reduced subvariety defined
by the radical ideal
    $\frq$. The symbolic powers
    \[
        \frq^{(k)} \ \defeq \ \set{f \in k[X] \,|\,
\ord_z f \geq k,\ z\in Z \text{ generic }}
    \]
    form a graded system.\footnote{When $Z$ is
reducible, we ask that the condition hold at a general point of each component.
The fact that this is equivalent to the usual algebraic definition is a theorem
of Zariski and Nagata: see \cite[Chapter 3]{Eisenbud}.} \item[(iii).] Let
    $<$ be a term order on $k[x_1,\ldots,x_n]$ and
$\bb$ be an ideal. Then
    \[
        \aa_k \ \defeq \ \op{in}_<(\bb^k)
    \]
    defines a graded system of monomial ideals, where $\op{in}_<(\bb^k)$
    denotes the initial ideal with respect to the given term order.

\end{enumerate}
\end{example}

\begin{example} [Valuation ideals]
\label{Valuation.GSI.Ex} Let     $\nu$ be a
$\RR$-valued valuation centered on
$k[X]$. Then the valuation ideals
    \[
        \aa_k  \ \defeq\  \set{f \in k[X] \,|\,
\nu(f)\geq k}
    \]
    form a graded family. Special cases of this
construction are interesting even when $X = \AA^2_\CC$.
\begin{enumerate}
   \item[(i).]
    Let $\eta: Y \to \AA^2$ be a birational map with
$Y$ also smooth and let
$E
    \subseteq Y$ be a prime divisor. Define the
valuation $\nu(f) \defeq
    \ord_E(f)$. Then
    \[
        \aa_k \defeq \mu_* \Oo_Y(-kE) = \set{f \in
\Oo_X\,|\,
\nu(f)=\ord_E(f) \geq k}.
    \]
     \item[(ii).]
    In $\CC[x,y]$ put $\nu(x)=1$ and
    $\nu(y)=\frac{1}{\sqrt{2}}$. Then one gets a
valuation by weighted degree. Here
$\aa_k$
    is the monomial ideal generated by the monomials
$x^iy^j$ such that
$i+\tfrac{j}{\sqrt{2}}
    \geq k$.
    \item[(iii).]
    Given $f \in \CC[x,y]$ define
$\nu(f)=\ord_z(f(z,e^z-1))$. This
    yields a valuation giving rise to the graded
system
    \[
        \aa_k \ \defeq \ (x^k,y-P_{k-1}(x))
    \]
    where $P_{k-1}(x)$ is the $(k-1)$st Taylor polynomial
of $e^x-1$. Note that the general element in $\fra_k$
defines a smooth curve in the plane.
\end{enumerate}
\end{example}

\begin{remark}
Except for Example \ref{First.Ex.GSI}(i), all these
constructions give graded families $\fra_{\bullet}$ for
which the corresponding Rees algebra need not be
finitely generated.
\end{remark}

\subsection{Asymptotic multiplier ideals} We now attach
multiplier ideals $\MI{\fra\bull^c}$ to a graded
family
$\fra\bull$ of ideals. The starting point is:
\begin{lemma} \label{Def.AMI.Lemma}
    Let $\aa_\bullet$ be  a graded system of ideals on
    $X$, and fix a rational number $c > 0$. Then for $p
\gg 0$  the multiplier ideals
    $\Jj(\aa^{c/p}_{ p})$ all coincide.
\end{lemma}

\begin{definition} \label{Def.AMI}
    Let $\aa_\bullet=\set{\aa_k}_{k \in \NN}$ be a
graded system of ideals on
    $X$. Given $c > 0$ we define the
\emph{asymptotic multiplier ideal} of $\fra\bull$ with
exponent $c$ to be the common ideal
    \[
          \Jj(\aa_\bullet^c) \ \defeq \
\Jj(\aa^{c/p}_{ p})
    \]
    for any sufficiently big $p\gg 0$.\footnote{In
\cite{ELS1} and early versions of \cite{PAG}, one only
dealt with the ideals $\MI{\fra\bull^{\ell}} $ for
integral $\ell$, which were written
$\MI{\Vert \fra_\ell \Vert}$.}
\end{definition}

\begin{proof}[Indication of Proof of Lemma
\ref{Def.AMI.Lemma}]
    We first claim that one has an inclusion of
multiplier
    ideals $\Jj(\aa^{c/p}_{ p}) \subseteq
\Jj(\aa^{c/pq}_{ pq})$ for all
    $p,q \geq 0$. Granting this, it follows from the
Noetherian condition that the collection of ideals
$\set{\Jj(\aa^{c/p}_{p})}_{p\geq 0}$
    has a unique maximal element. This proves
    the lemma at least for sufficiently divisible
$p$. (The statement for all $p \gg 0$ requires a little more
work; see \cite[Chapter 11]{PAG}.)

    To verify the claim let $\mu: X^\prime \lra X$ be a
common log resolution of
    $\aa_{ p }$ and $\aa_{  pq}$ with
$\aa_{  p}
\cdot
\Oo_Y =
\Oo_Y(-F_{  p})$ and
    $\aa_{  pq} \cdot \Oo_Y = \Oo_Y(-F_{  pq})$.
Since the $\fra_k$ form a graded
    system one has $\aa_{  p}^q \subseteq
\aa_{ pq}$ and therefore
    $-cqF_{  p}\leq-cF_{  pq}$. Thus
    \[
\mu_*\Oo_Y(K_{Y/X}-\floor{\textstyle{\frac{cq}{pq}}F_
{  p}})
\ \subseteq \
\mu_*\Oo_Y(K_{Y/X}-\floor{\frac{c}{pq}F_{ pq}})
    \]
    as claimed.
\end{proof}

\begin{remark}
Lemma \ref{Def.AMI.Lemma} shows that any information
captured by the multiplier ideals $\MI{\fra_{
p}^{c/p}}$ is present already for any one sufficiently
large index $p$. It is in this sense that multiplier
ideals have some finiteness built in that may not be
present in the underlying graded system $\fra\bull$.
\end{remark}

\begin{exercise} \label{Ex.AMI.computations}
    We return to the graded systems in Example
\ref{Valuation.GSI.Ex} coming from  valuations on
$\AA^2$.
\begin{enumerate}
\item [(ii).] Here $\fra_k $ is the monomial ideal
generated by $x^i y^j$ with $i + \tfrac{j}{\sqrt 2} \ge
k$, and
    $\Jj(\fra\bull^c)$ is the
monomial ideal generated by all  $x^iy^j$ with
\[ (i+1) +
\frac{(j+1)}{\sqrt{2}} \ > \ c .\]
(Compare with Theorem \ref{thm.mon}.)
\item[(iii).]
    Now take the valuation $\nu(f) = \ord_z
f(z,e^z-1)$. Then
    \[
        \Jj(\fra\bull^c) \ = \ \CC[x,y]
    \]
    for all $c > 0$. (Use the fact that each
$\aa_k$ contains a smooth curve.)
\end{enumerate}
\end{exercise}

\subsection{Growth of graded systems}
We now use the Subadditivity Theorem
\ref{thm.subadd} to prove a result from
\cite{ELS1} concerning the multiplicative behavior of
graded families of ideals:
\begin{theorem}\label{thm.asympIncl}
    Let $\aa_\bullet$ be a graded system of ideals and
fix any $\ell \in \NN$. Then
\[ \MI{\fra\bull^\ell} \ = \ \MI{\fra_{\ell
p}^{1/p}} \ \ \text{for}  \ \ p \gg 0.\]
Moreover for every
$m \in
\NN$ one has:
    \begin{equation} \label{Sub.AMI.Eqn}
        \aa_\ell^m \ \subseteq \ \aa_{\ell m} \subseteq
\ \Jj(\fra\bull^{\ell m}) \
\subseteq \
        \Jj(\fra\bull^\ell)^m.
    \end{equation}
    In particular, if $\Jj(\fra\bull^\ell) \subseteq
\bb$ for some natural number
$\ell$ and
    ideal $\bb$, then $\aa_{\ell m} \subseteq \bb^m$ for
all
$m$.
\end{theorem}

\begin{remark}
The crucial point here is the containment $\Jj(\fra\bull^{\ell m}) \subseteq
        \Jj(\fra\bull^\ell)^m$: it shows that passing
to multiplier ideals ``reverses" the inclusion $\fra_\ell^m \subseteq
\fra_{\ell m}$.
\end{remark}
\begin{proof} [Proof of Theorem \ref{thm.asympIncl}]
For the first statement, observe  that if $p \gg 0$ then
\[ \MI{\fra\bull^\ell} \ = \ \MI{\fra_p^{\ell/p}} \ = \
\MI{\fra_{ \ell p}^{\ell/ \ell p}} \ = \ \MI{\fra^{1/p}_{\ell p}}, \] where the
second equality is obtained by taking $  \ell p$ in place of  $p$ as the large
index in Lemma \ref{Def.AMI.Lemma}. For the containment $ \fra_{\ell m}
\subseteq \MI{\fra\bull^{\ell m}}$ it is then enough to prove that $\fra_{\ell
m} \subseteq \MI{\fra_{\ell m p}^{1/p}}.$ But we have $\fra_{\ell m} \subseteq
\MI{\fra_{\ell m}}$ thanks to Exercise \ref{ex.intclo}, while the inclusion
$\MI{\fra_{\ell m}} \subseteq \MI{\fra^{1/p}_{\ell m p}}$ was established
during the proof of \ref{Def.AMI.Lemma}.

It remains only to prove that $\MI{\fra\bull^{\ell m}}
\subseteq \MI{\fra\bull^\ell}^m$.  To this end, fix
$p
\gg 0$. Then by the definition of asymptotic multiplier
ideals   and the Subadditivity Theorem one has
\begin{align*} \MI{\fra\bull^{\ell m}} \ &= \
          \Jj(\aa^{\ell m/p}_{p}) \\  &\subseteq \
        \MI{\fra_{ p}^{\ell/p}}^m \\
&= \ \MI{\fra\bull^\ell}^m,
    \end{align*}
as required.\end{proof}

\begin{example}
The Theorem gives another explanation of the fact  that the multiplier ideals
associated to the graded system $\fra\bull$ from Example
\ref{Valuation.GSI.Ex}.(iii) are trivial. In fact, in this example the colength
of $\fra_k$ in $\CC[X]$ grows linearly in $k$. It follows from Theorem
\ref{thm.asympIncl} that then $\MI{\fra\bull^\ell} = (1)$ for all $\ell$.
\end{example}

\begin{exercise}\label{Triv.GSI.AMI}
 Let $\fra_k = \frb^k$ be
the trivial graded family consisting of powers of a
fixed ideal. Then $\MI{\fra\bull^c} = \MI{\frb^c}$ for
all $c > 0$. So we do not get anything new in this
case.
\end{exercise}

\subsection{A comparison theorem for symbolic powers}\label{sec.symb}

As a quick but surprising application of Theorem \ref{thm.asympIncl} we discuss
a result due to Ein, Smith and the second author from \cite{ELS1}
concerning symbolic powers of radical ideals.

Consider  a reduced subvariety $Z \subseteq X$ defined by a radical ideal $\frq
\subseteq k[X]$. Recall from Example \ref{First.Ex.GSI}(ii) that one can define
the symbolic powers $\frq^{(k)}$ of $\frq$ to be
   \[
        \frq^{(k)} \ \defeq \ \set{f \in \Oo_X \,|\,
\ord_z f \geq k,\ z\in Z}.
    \]
Thus evidently $\frq^k \subseteq \frq^{(k)}$, and
equality holds if $Z$ is smooth. However if $Z$ is
singular then in  general the inclusion is strict:
\begin{example}
Take $Z \subseteq \CC^3$ to be the
union of the three coordinate axes, defined by the
ideal
\[ \frq \ = \ ( \, xy \, , \, yz \, , \, xz \,) \
\subseteq \ \CC[x,y,z]. \] Then  $xyz \, \in \, \frq^{(2)}$ since evidently the
union of the three coordinate planes has multiplicity $2$ at a general point of
$Z$. But $\frq^2$ is generated by monomials of degree $4$, thus cannot contain
$xyz$, which is of degree $3$.
\end{example}

Swanson \cite{Swanson00a} proved (in a much more general setting)
that there exists an integer $k = k(Z)$ such that
\[ \frq^{(km)} \ \subseteq \ \frq^m \]
for all $m \ge 0$. At first glance, one might be
tempted to suppose
that for very singular $Z$ the coefficient $k(Z)$ will
have to become quite large. The main result of
\cite{ELS1} shows that this isn't the case, and that in
fact one can take
$k(Z) =
\codim Z$:
\begin{theorem} \label{Symb.Pow.Thm} Assume that every
irreducible component of $Z$ has codimension $\le e$ in
$X$. Then
\[ \frq^{(e m)} \ \subseteq \frq^m \ \text{ for all } \
m \ge 0. \]
In particular, $\frq^{(m \cdot \dim X)} \subseteq
\frq^m$ for all radical ideals $\frq \subseteq k[X]$
and all $m
\ge 0$.
\end{theorem}

\begin{example} [Points in the plane] Let $T \subseteq
\PP^2$ be a finite set (considered as a reduced scheme), and let $I \subseteq S
= \CC[x,y,z]$ be the homogeneous ideal of $T$. Suppose that $f \in S$ is a
homogeneous form which has multiplicity $\ge 2m$ at each of the points of $T$.
Then $f \in I^m$. (Apply Theorem \ref{Symb.Pow.Thm} to the homogeneous ideal
$I$ of $T$.) In spite of the classical nature of this statement, we do not know
a direct elementary proof.
\end{example}

\begin{proof} [Proof of Theorem \ref{Symb.Pow.Thm}]
Applying Theorem \ref{thm.asympIncl} to the graded
system $\fra_k =
\frq^{(k)}$, it suffices to show that
\begin{equation}
\MI{\fra\bull^e} \ \subseteq \ \frq. \tag{*}
\end{equation}
Since $\frq$ is radical, it suffices to test the
inclusion (*) at a general point of $Z$. Therefore we
can assume that $Z$ is smooth, in which case
$\frq^{(k)} = \frq^k$. Now Exercises
\ref{Smooth.Ideal.Exercise}   and \ref{Triv.GSI.AMI}
apply.
\end{proof}

\begin{remark}
Using their theory of tight closure,
Hochster and Huneke \cite{HH.ComOrdSymbPow} have extended Theorem
\ref{Symb.Pow.Thm} to arbitrary regular Noetherian rings containing a
field.
\end{remark}

\begin{remark}
Theorem \ref{thm.asympIncl}  is applied in \cite{ELS2} to study the
multiplicative behavior of Abyhankar valuations centered at a smooth point of a
complex variety.
\end{remark}

\begin{remark}
Working with the asymptotic multiplier ideals $\MI{\fra\bull^c}$ one can define
the log-canonical threshold and jumping coefficients of a graded system
$\fra\bull$ much as in \S 4. However now these numbers need no longer be
rational, the periodicy of jumping numbers (Exercise
\ref{Periodicityu.Jump.Nos}) may fail, and in fact the collection of jumping
coefficients of $\fra\bull$ can contain accumulation points. See \cite[\S
5]{ELSV}.
\end{remark}

\bibliographystyle{amsalpha}
\bibliography{MultiplierNotes}

\providecommand{\bysame}{\leavevmode\hbox to3em{\hrulefill}\thinspace}
\providecommand{\MR}{\relax\ifhmode\unskip\space\fi MR }
\providecommand{\MRhref}[2]{%
  \href{http://www.ams.org/mathscinet-getitem?mr=#1}{#2}
}
\providecommand{\href}[2]{#2}
\begin{thebibliography}{DEL00}

\bibitem[AS95]{Angehrn-Siu95a}
Urban Angehrn and Yum~Tong Siu, \emph{Effective freeness and point separation
  for adjoint bundles}, Invent. Math. \textbf{122} (1995), no.~2, 291--308.

\bibitem[Bli]{Bli.MultToric}
Manuel Blickle, \emph{{Multiplier ideal and module on a toric variety.}},
  preprint.

\bibitem[BM97]{Bierstone-Milman97}
Edward Bierstone and Pierre~D. Milman, \emph{Canonical desingularization in
  characteristic zero by blowing up the maximum strata of a local invariant},
  Invent. Math. \textbf{128} (1997), no.~2, 207--302.

\bibitem[Dan78]{Danilov78}
V.~I. Danilov, \emph{The geometry of toric varieties}, Russian Math. Surveys
  \textbf{33} (1978), 97--154.

\bibitem[DEL00]{DEL}
Jean-Pierre Demailly, Lawrence Ein, and Robert Lazarsfeld, \emph{A
  subadditivity property of multiplier ideals}, Michigan Math. J. \textbf{48}
  (2000), 137--156, Dedicated to William Fulton on the occasion of his 60th
  birthday.

\bibitem[Dem93]{Demailly93c}
Jean-Pierre Demailly, \emph{A numerical criterion for very ample line bundles},
  J. Differential Geom. \textbf{37} (1993), no.~2, 323--374.

\bibitem[Dem99]{Demailly99b}
\bysame, \emph{M\'ethodes ${L}^2$ et r\'esultats effectifs en g\'eom\'etrie
  alg\'ebrique}, S\'eminaire Bourbaki (1999).

\bibitem[Dem01]{Dem.Mult}
\bysame, \emph{Multiplier ideal sheaves and analytic methods in algebraic
  geometry}, School on Vanishing Theorems and Effective Results in Algebraic
  Geometry (Trieste, 2000), ICTP Lect. Notes, vol.~6, Abdus Salam Int. Cent.
  Theoret. Phys., Trieste, 2001, pp.~1--148. \MR{1 919 457}

\bibitem[Eis95]{Eisenbud}
David Eisenbud, \emph{{Commutative algebra. With a view toward algebraic
  geometry.}}, {Graduate Texts in Mathematics. 150. Berlin: Springer-Verlag.
  xvi, 785 p.}, 1995.

\bibitem[EL97]{Ein-Lazarsfeld97a}
Lawrence Ein and Robert Lazarsfeld, \emph{Singularities of theta divisors and
  the birational geometry of irregular varieties}, J. Amer. Math. Soc.
  \textbf{10} (1997), no.~1, 243--258.

\bibitem[EL99]{ELNull}
\bysame, \emph{A geometric effective {N}ullstellensatz}, Invent. Math.
  \textbf{137} (1999), no.~2, 427--448.

\bibitem[ELS]{ELS2}
Lawrence Ein, Robert Lazarsfeld, and Karen~E. Smith, \emph{Uniform
  approximation of abhyankar valuations on smooth function fields}, Am. J.
  Math, to appear.

\bibitem[ELS01]{ELS1}
\bysame, \emph{Uniform bounds and symbolic powers on smooth varieties}, Invent.
  Math. \textbf{144} (2001), no.~2, 241--252.

\bibitem[ELSV]{ELSV}
Lawrence Ein, Robert Lazarsfeld, K.~E. Smith, and Dror Varolin, \emph{Jumping
  coefficients of multiplier ideals}, To appear.

\bibitem[EV92]{EV}
H{\'e}l{\`e}ne Esnault and Eckart Viehweg, \emph{Lectures on vanishing
  theorems}, DMV Seminar, vol.~20, Birkh\"auser Verlag, Basel, 1992.

\bibitem[EV00]{EncVill.Desing}
Santiago Encinas and Orlando Villamayor, \emph{A course on constructive
  desingularization and equivariance}, Resolution of singularities (Obergurgl,
  1997), Progr. Math., vol. 181, Birkh\"auser, Basel, 2000, pp.~147--227.

\bibitem[Ful93]{Fulton.Toric}
William Fulton, \emph{Introduction to toric varieties}, Annals of Mathematics
  Studies, vol. 131, Princeton University Press, Princeton, NJ, 1993, The
  William H. Roever Lectures in Geometry.

\bibitem[Gow02]{Goward.PrincMon}
Russel A.~Jr. Goward, \emph{A simple algorithm for principalization of monomial
  ideals}, preprint, 2002.

\bibitem[Har77]{Ha}
Robin Hartshorne, \emph{Algebraic geometry}, Springer-Verlag, New York, 1977,
  Graduate Texts in Mathematics, No. 52.

\bibitem[Har01]{Hara.GeomIntTest}
Nobuo Hara, \emph{{Geometric interpretation of tight closure and test
  ideals.}}, Trans. Am. Math. Soc. \textbf{353} (2001), no.~5, 1885--1906
  (English).

\bibitem[HH90]{HH.BrianSkoda}
Melvin Hochster and Craig Huneke, \emph{Tight closure, invariant theory, and
  the {B}rian\c con-{S}koda theorem}, J. Amer. Math. Soc. \textbf{3} (1990),
  no.~1, 31--116.

\bibitem[HH02]{HH.ComOrdSymbPow}
\bysame, \emph{Comparison of symbolic and ordinary powers of ideals}, Invent.
  Math. \textbf{147} (2002), no.~2, 349--369.

\bibitem[Hir64]{Hironaka.ResSing}
Heisuke Hironaka, \emph{Resolution of singularities of an algebraic variety
  over a field of characteristic zero. {I}, {II}}, Ann. of Math. (2) 79 (1964),
  109--203; ibid. (2) \textbf{79} (1964), 205--326.

\bibitem[How01]{Howald.MultMon}
J.~A. Howald, \emph{Multiplier ideals of monomial ideals}, Trans. Amer. Math.
  Soc. \textbf{353} (2001), no.~7, 2665--2671 (electronic).

\bibitem[HT]{HaTak}
Nobuo Hara and Shunsuke Takagi, \emph{{Some remarks on a generalization of test
  ideals}}.

\bibitem[Hun92]{Huneke.UniformBounds}
Craig Huneke, \emph{Uniform bounds in {N}oetherian rings}, Invent. Math.
  \textbf{107} (1992), no.~1, 203--223.

\bibitem[HY]{HaraYosh}
{Nobuo} {Hara} and {Ken-ichi} Yoshida, \emph{{A generalization of tight closure
  and multiplier ideals}}, arXiv:math.AC/0211008.

\bibitem[Kol97]{Kollar.Sings.of.Pairs}
J{\'a}nos Koll{\'a}r, \emph{Singularities of pairs}, Algebraic geometry---Santa
  Cruz 1995, Proc. Sympos. Pure Math., vol.~62, Amer. Math. Soc., Providence,
  RI, 1997, pp.~221--287.

\bibitem[Laz]{PAG}
Robert Lazarsfeld, \emph{Positivity in algebraic geometry}, Monograph in
  preparation.

\bibitem[Lib83]{Libgober}
A.~Libgober, \emph{Alexander invariants of plane algebraic curves},
  Singularities, Part 2 (Arcata, Calif., 1981), Proc. Sympos. Pure Math.,
  vol.~40, Amer. Math. Soc., Providence, RI, 1983, pp.~135--143.

\bibitem[Lip93]{lip.adj}
Joseph Lipman, \emph{Adjoints and polars of simple complete ideals in
  two-dimensional regular local rings}, Bull. Soc. Math. Belg. S\'er. A
  \textbf{45} (1993), no.~1-2, 223--244, Third Week on Algebra and Algebraic
  Geometry (SAGA III) (Puerto de la Cruz, 1992).

\bibitem[LV90]{Loeser.Vaquie}
F.~Loeser and M.~Vaqui{\'e}, \emph{Le polyn\^ome d'{A}lexander d'une courbe
  plane projective}, Topology \textbf{29} (1990), no.~2, 163--173.

\bibitem[Nad90]{Nadel90}
Alan~Michael Nadel, \emph{Multiplier ideal sheaves and {K}\"ahler-{E}instein
  metrics of positive scalar curvature}, Ann. of Math. (2) \textbf{132} (1990),
  no.~3, 549--596.

\bibitem[Par99]{Paranjape}
Kapil~H. Paranjape, \emph{The {B}ogomolov-{P}antev resolution, an expository
  account}, New Trends in Algebraic Geometry (Warwick, 1996), London Math. Soc.
  Lecture Note Ser., vol. 264, Cambridge Univ. Press, Cambridge, 1999,
  pp.~347--358. \MR{1 714 830}

\bibitem[SB74]{BS.74}
Henri Skoda and Jo{\"e}l Brian{\c{c}}on, \emph{Sur la cl\^oture int\'egrale
  d'un id\'eal de germes de fonctions holomorphes en un point de ${\bf
  {c}}\sp{n}$}, C. R. Acad. Sci. Paris S\'er. A \textbf{278} (1974), 949--951.

\bibitem[Siu98]{Siu98a}
Yum-Tong Siu, \emph{Invariance of plurigenera}, Invent. Math. \textbf{134}
  (1998), no.~3, 661--673.

\bibitem[Sko72]{Skoda72}
Henri Skoda, \emph{Application des techniques ${L}\sp{2}$ \`a la th\'eorie des
  id\'eaux d'une alg\`ebre de fonctions holomorphes avec poids}, Ann. Sci.
  \'Ecole Norm. Sup. (4) \textbf{5} (1972), 545--579.

\bibitem[Smi00]{Smith.MultTest}
Karen~E. Smith, \emph{{The multiplier ideal is a universal test ideal.}},
  Commun. Algebra \textbf{28} (2000), no.~12, 5915--5929 (English).

\bibitem[Swa00]{Swanson00a}
Irena Swanson, \emph{Linear equivalence of ideal topologies}, Math. Z.
  \textbf{234} (2000), no.~4, 755--775.

\bibitem[Taka]{Takagi.MultTest}
Shunsuke Takagi, \emph{{An interpretation of multiplier ideals via tight
  closure}}, arXiv:math.AG/0111187.

\bibitem[Takb]{Tak.Inv}
\bysame, \emph{{F-singularities of pairs and Inversion of Adjunction of
  arbitrary codimension}}.

\bibitem[TiW]{TakWat}
Shunsuke Takagi and Kei ichi Watanabe, \emph{{When does the subadditivity
  theorem for multiplier ideals hold?}}

\end{thebibliography}
\end{document}